\documentclass[12pt,leqno,draft]{article}
\usepackage{amsfonts}
\usepackage{amsmath, amsthm, amsfonts, amssymb, color}
\usepackage{mathrsfs}
\usepackage{color}
\usepackage[notcite,notref]{showkeys}
\usepackage{cancel}

\newtheorem{thm}{Theorem}[section]
\newtheorem{cor}[thm]{Corollary}

\theoremstyle{definition}

\newcommand{\scr}[1]{\mathscr #1}
\definecolor{wco}{rgb}{0.5,0.2,0.3}

\numberwithin{equation}{section} \theoremstyle{remark}

\newcommand{\ua}{\uparrow}

\title{{\bf Regularity Estimates  for Singular Density Dependent SDEs}
	\footnote{Supported in part by the National Key R\&D Program of China (2022YFA1006000) and NSFC(12531007, 12101390, 12426656).} }
\author{
	{\bf    Feng-Yu Wang$^{(a)}$, Qiumiao Wen$^{(b,c)}$ ,  Fen-Fen Yang$^{(b,c)}$   }\\
	\footnotesize{
		(a)  Center for Applied Mathematics and KL-AAGDM, Tianjin
		University, Tianjin 300072, China }\\
	\footnotesize{ (b) 	Department of Mathematics, Shanghai University, Shanghai 200444, China  }\\
	\footnotesize{ (c) Newtouch Center for Mathematics, Shanghai University, Shanghai, 200444, China}\\
	\footnotesize{  wangfy@tju.edu.cn; wen\_qm@shu.edu.cn;   yangfenfen@shu.edu.cn }
}
\begin{document}
\allowdisplaybreaks
\def\R{\mathbb R}  \def\ff{\frac} \def\ss{\sqrt} \def\B{\mathbf
B} \def\W{\mathbb W}
\def\N{\mathbb N} \def\kk{\kappa} \def\m{{\bf m}}
\def\ee{\varepsilon}\def\ddd{D^*}
\def\dd{\delta} \def\DD{\Delta} \def\vv{\varepsilon} \def\rr_V{\rho}
\def\<{\langle} \def\>{\rangle} \def\GG{\Gamma} \def\gg{\gamma}
  \def\nn{\nabla} \def\pp{\partial} \def\E{\mathbb E}
\def\d{\text{\rm{d}}} \def\bb{\beta} \def\aa{\alpha} \def\D{\scr D}
  \def\si{\sigma} \def\ess{\text{\rm{ess}}}
\def\beg{\begin} \def\beq{\begin{equation}}  \def\F{\scr F}
\def\Ric{\text{\rm{Ric}}} \def\Hess{\text{\rm{Hess}}}
\def\e{\text{\rm{e}}} \def\ua{\underline a} \def\OO{\Omega}  \def\oo{\omega}
 \def\tt{\tilde} \def\Ric{\text{\rm{Ric}}}
\def\cut{\text{\rm{cut}}} \def\P{\mathbb P} \def\ifn{I_n(f^{\bigotimes n})}
\def\C{\scr C}      \def\aaa{\mathbf{r}}     \def\r{r}
\def\gap{\text{\rm{gap}}} \def\prr{\pi_{{\bf m},\varrho}}  \def\r{\mathbf r}
\def\Z{\mathbb Z} \def\vrr{\varrho} \def\ll{\lambda}
\def\L{\scr L}\def\Tt{\tt} \def\TT{\tt}\def\II{\mathbb I}
\def\i{{\rm in}}\def\Sect{{\rm Sect}}  \def\H{\mathbb H}
\def\M{\scr M}\def\Q{\mathbb Q} \def\texto{\text{o}} \def\LL{\Lambda}
\def\Rank{{\rm Rank}} \def\B{\scr B} \def\i{{\rm i}} \def\HR{\hat{\R}^d}
\def\to{\rightarrow}\def\l{\ell}\def\iint{\int}
\def\EE{\scr E}\def\Cut{{\rm Cut}}
\def\A{\scr A} \def\Lip{{\rm Lip}}
\def\BB{\scr B}\def\Ent{{\rm Ent}}\def\L{\scr L}
\def\R{\mathbb R}  \def\ff{\frac} \def\ss{\sqrt} \def\B{\mathbf
B}
\def\N{\mathbb N} \def\kk{\kappa} \def\m{{\bf m}}
\def\dd{\delta} \def\DD{\Delta} \def\vv{\varepsilon} \def\rr{\rho}
\def\<{\langle} \def\>{\rangle} \def\GG{\Gamma} \def\gg{\gamma}
  \def\nn{\nabla} \def\pp{\partial} \def\E{\mathbb E}
\def\d{\text{\rm{d}}} \def\bb{\beta} \def\aa{\alpha} \def\D{\scr D}
  \def\si{\sigma} \def\ess{\text{\rm{ess}}}
\def\beg{\begin} \def\beq{\begin{equation}}  \def\F{\scr F}
\def\Ric{\text{\rm{Ric}}} \def\Hess{\text{\rm{Hess}}}
\def\e{\text{\rm{e}}} \def\ua{\underline a} \def\OO{\Omega}  \def\oo{\omega}
 \def\tt{\tilde} \def\Ric{\text{\rm{Ric}}}
\def\cut{\text{\rm{cut}}} \def\P{\mathbb P} \def\ifn{I_n(f^{\bigotimes n})}
\def\C{\scr C}      \def\aaa{\mathbf{r}}     \def\r{r}
\def\gap{\text{\rm{gap}}} \def\prr{\pi_{{\bf m},\varrho}}  \def\r{\mathbf r}
\def\Z{\mathbb Z} \def\vrr{\varrho} \def\ll{\lambda}
\def\L{\scr L}\def\Tt{\tt} \def\TT{\tt}\def\II{\mathbb I}
\def\i{{\rm in}}\def\Sect{{\rm Sect}}  \def\H{\mathbb H}
\def\M{\scr M}\def\Q{\mathbb Q} \def\texto{\text{o}} \def\LL{\Lambda}
\def\Rank{{\rm Rank}} \def\B{\scr B} \def\i{{\rm i}} \def\HR{\hat{\R}^d}
\def\to{\rightarrow}\def\l{\ell}\def\BB{\mathbb B}
\def\8{\infty}\def\I{1}\def\U{\scr U} \def\n{{\mathbf n}}\def\v{V}\def\LL{{\bf L}}
\def\les{\lesssim}
\def\cor{\textcolor{red}}
\maketitle

\begin{abstract}  Consider the density dependent (i.e. Nemytskii-type) SDEs on $\R^d$, where  the drift $b_t(x,\rr(x),\rr)$
 is locally integrable in $(t,x)\in [0,\infty)\times\R^d$ and  may be singular in  the distribution density function $\rr$.
The  relative/Renyi entropies  between two   time-marginal distributions are estimated  by  using the Wasserstein distance  of initial distributions.
   When $d=1$ and $b_t$ decays at $t=0$ with rate $t^{\ff 1 2+}$,  our the relative entropy estimate  coincides with  the classical entropy-cost inequality for elliptic diffusion processes.
To estimate the Renyi entropy, a refined Khasminskii estimate is presented for singular SDEs which may be interesting by itself.

 \end{abstract} \noindent
\noindent
 Keywords: Density dependent SDEs, relative entropy, Renyi entropy, Wasserstein distance, $\tt L^k$-distance.

\section{Introduction}

Let $\scr P$ be the set of all probability measures on $\R^d$ equipped with the weak topology,  let $\ell_\xi$ be the distribution density function of a random variable with respect to the Lebesgue measure,
let $\B(\R^d)$ be the space of Borel measurable functions on $\R^d$, and let
$$\D_1^+:= \bigg\{\rr\in \B(\R^d):\ \rr\ge 0,\ \int_{\R^d} \rr(x)\d x= 1\bigg\}$$
which is a Polish space under the $L^1$-distance  $\|f-g\|_{1}:= \int_{\R^d}|f-g|(x)\d x.$

Consider the  following Nemytskii-type density dependent SDE on $\R^d$:
\beq\label{E0} \d X_t= b_t(X_t,\ell_{X_t}(X_t),\ell_{X_t})\d t+\si_t(X_t)\d W_t,\ \ \ t\in [0,T],\end{equation}
where $T\in (0,\infty)$ is a fixed time, $(W_t)_{t\in [0,T]}$ is an $m$-dimensional Brownian motion on a probability base (i.e. complete filtered probability space) $(\OO, \F, \{\F_t\}_{t\in [0,T]},\P),$
  and
$$b: (0,T]\times \R^d\times [0,\infty)\times \D_1^+ \to \R^d,\ \ \ \si: (0,T]\times\R^d\to \R^{d\otimes m}$$
are measurable.   We take $(0,T]$ for time interval to allow singular distributions of $X_0$.

According to Kac's propagation of chaos, the distribution dependence in a stochastic equation refers  to the mean-field interactions in the associated particle systems. A typical  interaction
 is of type
$$b(x,\mu):= \int_{\R^d} V(x-y)\mu(\d y),\ \ \ \mu\in \hat{\scr P},\ x\in \R^d,$$
 where the    interaction kernel  $V: \R^d\to\R^d$ is measurable and  $\hat{\scr P}$ is a subclass of $\scr P$ such that the integral exists. In particular, the   Coulomb kernel
$$V(x):= \ff{cx}{|x|^d},\ \ x\ne 0$$
  for some constant $c\in (0,\infty)$ describes electrostatic interactions between  charged particles.

  The interactions included in the density dependent SDE \eqref{E0} are  given by the Dirac function $\delta$, since
  the distribution density function $\rr_\mu$ for an absolutely continuous probability measure $\mu$  satisfies
   $$\rr_\mu(x)=  \int_{\R^d} \dd(x-y)\rr_\mu(y)\d y=\int_{\R^d} {\bf \dd}(x-y) \mu(\d y).$$
   In physics, the Dirac function  characterizes densities of particles and point charges. Since the Dirac function  beyonds the class of real functions,  it is more singular than the above mentioned interaction kernels.
     Comparing with the existing  literature of McKean-Vlasov SDEs depending on global distribution properties (e.g. expectations),
     much less is known on the density dependent SDE  \eqref{E0}.
 See \cite{BR1, BR2, HR,IR,JM, W23, Z} and references therein for the study of well-posedness, superposition principle, propagation of chaos and the Euler scheme. However, crucial properties like gardient/entropy estimates and ergodicity are not yet studied  for such a singular model.

 In recent years, the entropy-cost  inequality (also called log-Harnack inequality)
 $$\Ent(P_t^*\mu|P_t^*\nu)\le \ff c t \W_2(\mu,\nu)^2$$
 has been established for McKean-Vlasov SDEs,   see for instance \cite{RW25, HRW25} and references therein. Here,  $P_t^*\mu$ is the distribution of the solution at time $t$ with initial distribution $\mu$,
 $\W_2$ is the $2$-Wasserstein distance, and
  $\Ent$ is the relative entropy; namely, for two probability measures $\mu$ and $\nu$,
 $$\Ent(\mu|\nu):= \beg{cases} \int_{\R^d}(\log \ff{\d\mu}{\d\nu})\d\mu, &\text{if}\  \ff{\d\mu}{\d\nu} \ \text{exists},\\
 \infty, &\text{otherwise.}\end{cases}$$
  On the other hand, the Renyi entropy
  $$\Ent_\aa(\mu|\nu):= \ff 1\aa \log\int_{\R^d}\Big(\ff{\d\mu}{\d\nu}\Big)^\aa \d\mu,\ \ \ \aa>0$$
  for   $P_t^*$ has been characterized by the dimension-free Harnack inequality with power developed by the first named  author, see \cite[Theorems 1.4.1 and Theorem 1.4.2] {W13} for  applications of log-Harnack and dimension-free Harnack inequalities.

 In this paper, we estimate the relative entropy  and Renyi entropy for time-marginal distributions for solutions of the density dependent SDE \eqref{E0}, which is open so far.
 We will allow the drift $b_t(x,r,\rr)$ for $(t,x,r,\rr)\in [0,T]\times \R^d\times [0,\infty)\times \D_1^+$ to be singular in $(t,x)$ and $\rr$.

 To measure the singularity in $x\in\R^d$,
 let
$$\|f\|_{\tt L^k}:= \sup_{z\in \R^d} \|1_{B(z,1)} f\|_{k},\ \ k\ge 1$$ for a measurable function (or vector field) $f$ on $\R^d,$
where $B(z,1):=\{x\in\R^d: |x-z|\le 1\}$, and $\|\cdot\|_{k}$ is the $L^k$-norm with respect to the Lebesgue measure.

To measure the singularity in the distribution parameter, for $k\in (1,\infty]$ and   $\mu,\nu\in \scr P$ with densities $\rr_\mu$ and $\rr_\nu$, let
\beg{align*}  \|\mu\|_{\tt L^k}:= \|\rr_\mu\|_{\tt L^k},\ \ \ \|\mu-\nu\|_{\tt L^k}:= \|\rr_\mu-\rr_\nu\|_{\tt L^k},\ \ k\in (1,\infty];\end{align*}
while for $k=1$ and   $\mu,\nu\in \scr P$ let
\beg{align*} \|\mu\|_{\tt L^1}:=\sup_{z\in \R^d, |f|\le 1} |\mu(1_{B(z,1)}f)|,\ \ \ \|\mu-\nu\|_{\tt L^1}:= \sup_{z\in\R^d, |f|\le 1} |(\mu-\nu)(1_{B(z,1)}f)|.\end{align*}
Let
$$\tt{\scr P}_k:= \Big\{\mu\in \scr P:\ \|\mu\|_{\tt L^k}<\infty\Big\},\ \ k\in [1,\infty].$$
  In particular, when $k=1$ we have     $\tt{\scr P}_1=\scr P$ and
$$\|\mu-\nu\|_{\tt L^1}\le \|\mu-\nu\|_{var}:= \sup_{|f|\le 1} |\mu(f)-\nu(f)|.$$

For any  $p,q\in [1,\infty]$, let $\tt L_{q}^p(T)$ be the space of measurable functions $f$ on $[0,T]\times \R^d$ such that
$$\|f\|_{\tt L_q^p(T)}:= \sup_{z\in\R^d} \bigg(\int_0^T\|f_t1_{B(z,1)} \|_{p}^q\d t\bigg)^{\ff 1 q} <\infty.$$
We will take $(p,q)$ from the class
$$\scr K:=\Big\{(p,q)\in (2,\infty]^2:\ \ff d p+\ff 2 q<1\Big\}.$$

  Let $\L_{\xi}$ denote the distribution of a random variable $\xi$. If different probability measures are considered, we
write $\L_{\xi|\P}$ instead of $\L_{\xi}$ to emphasize the underlying probability $\P$. Let $X_{[0,T]}$ denote the stochastic process $[0,T]\ni t\mapsto X_t$.

 \beg{defn} Let $X_0$ be $\F_0$-measurable with $\mu:=\L_{X_0} \in\scr P.$
 \beg{enumerate}
 \item[$(1)$] We call $X_t$ a  (strong) solution of \eqref{E0} with initial value $X_0$, if it is a continuous adapted process with  such that $\ell_{X_t}$ exists for $t\in (0,T]$,
 and $\P$-a.s.
 $$X_t= X_0+\int_0^t b_s(X_s, \ell_{X_s}(X_s), \ell_{X_s})\d s+\int_0^t \si_s(X_s)\d W_s,\ \ t\in [0,T].$$
 \item[$(2)$] A pair $(\tt X_t, \tt W_t)$ is called a weak  solution of \eqref{E0} with initial distribution $\mu$, if there exists a    probability base $(\tt\OO, \tt\F, \{\tt\F_t\}_{t\in [0,T]},\tt \P)$ under which  $\tt W_t$ is an $m$-dimensional Brownian motion, $\L_{\tt X_0|\tt\P}=\mu$  and
 $\tt X_t$ is a  solution of \eqref{E0} with $(\tt X_t,\tt W_t)$ in place of $(X_t,W_t).$
  \item[$(3)$]  We say that \eqref{E0} has weak uniqueness  with
 initial distribution $\mu$, if for any two weak solutions $(X_t^i, W_t^i)$   under probability bases $(\OO^i, \F^i, \{\F_t^i\}_{t\in [0,T]},\P^i)$ with
 $\L_{X_0^i|\P^i}=\mu$ for $i=1,2,$
 we have $\L_{X_{[0,T]}^1|\P^1}=\L_{X_{[0,T]}^2|\P^2}.$
 \end{enumerate} \end{defn}

 If,   the SDE \eqref{E0} has a unique (weak)  solution  $X_t$ for  initial distribution $\mu=\L_{X_0}\in \tt {\scr P}_p$,  we denote
 $$ P_t^*\mu :=\L_{X_t},\ \ \ t\in [0,T].$$
%\cor {\beg{defn}  We call  the  solution  of \eqref{E0}  super-continuous in initial distribution, if
%	$$\lim_{\| \mu- \nu\|_{\tt L^p}\to 0} \|P_t^*\mu-P_t^*\nu\|_{\tt L^k}=0,\ \ \ t\in (0,T]  $$ holds for some constants $k>p\ge 1$.
%	\end{defn}}

In Section 2, we study the well-posedness of \eqref{E0} and estimate the super-continuity of $P_t^*$, i.e. for $p\in [1,k)$,
$$\|P_t^*\mu-P_t^*\nu\|_{\tt L^k}\le \aa_{p,k}(t) \|\mu-\nu\|_{\tt L^p},\ \ \ t>0,\ \mu,\nu\in \tt{\scr P_p}$$
holds for some rate function $\aa_{p,k}: (0,\infty)\to (0,\infty)$  with $\aa_{p,k}(t)\sim t^{-\ff{d(k-p)}{2kp}}$ for small $t>0$.
This estimate is then applied in Sections 3 and 4 to estimate the relative entropy and Renyi entropy between $P_t^*\mu$ and $P_t^*\nu$.

\section{Well-posedness and super-continuity}

 Let $a_t(x):= (\si_t\si_t^*)(x)$ and   decompose $b$ as
$$b_t(x,r, \rr)=   b_t^{(1)}(x)+\sum_{i=2}^l b_t^{(i)}(x,r,\rr),$$
 where $2\le l\in \mathbb N$ and $\{b^{(i)}\}_{1\le i\le l} $ are measurable in all arguments.
 %\cor{The solution of $\eqref{E0}$ is defined in a standard way?}, see \cite{HR,W23}.

\beg{enumerate}\item[{\bf (A)}]   \emph{
 There exist  $K\in (0,\infty), \vv\in (0,1),   2\le l\in \mathbb N,$    $l'\in\mathbb N$ and
 $$ (p_i,q_i),  \in  \scr K,\ \ 0\le f^{(i)},g^{(j)}\in \tt L_{q_i}^{p_i}(T), \ \ 2\le i\le l,\ \ 1\le j\le l',$$ such that   for all
$t\in (0,T], r\ge 0, \rr\in \D_1^+$ and $ x,y\in \R^d,$
\beq\label{P1} \|a \|_\infty+ \|a ^{-1}\|_\infty\le K,\ \  |a_t(x)- a_{t}(y)|\le K|x-y|^\vv,\  \ \|\nn a\|\le \sum_{j=1}^{l'} g^{(j)},
\end{equation}
\beq\label{P2}  \|\nn b_t^{(1)}\|_\infty
 +\|b_0^{(1)} \|_{\infty }  \le K,\ \ \qquad |b_t^{(i)} (x,r,\rr)|\le f_t^{(i)}(x),\ \ 2\le i\le l.\end{equation}
Moreover, there exists $k\in (d,\infty]$   such that for any $t\in (0,T],\ x\in \R^d,\ r,r'\in [0,\infty)$ and $\rr,\rr'\in \D_1^+$,
\beq\label{P3}   |b_t(x,r,\rr)-b_t(x,r',\rr')|\le K \big(|r-r'|+t^{-\ff{d}{2k}} \|\rr-\rr'\|_{\tt L^k}\big),\end{equation}}
\end{enumerate}

Recall that   the $q$-Wasserstein distance for $q\in [1,\infty)$ is defined as
$$\W_q(\mu,\nu):= \inf_{\pi\in \C(\mu,\nu)} \bigg(\int_{\R^d\times\R^d} |x-y|^q\pi(\d x,\d y)\bigg)^{\ff 1 q},$$
where $\C(\mu,\nu)$ is the set of all couplings of $\mu$ and $\nu$.

 \beg{thm}\label{T1} Assume {\bf (A)}.  Then the following assertions hold for any $p\in (d,k]\cap [1,\infty)$.
 \beg{enumerate} \item[$(1)$]
  For  any $\F_0$-measurable initial value $X_0$  with $\L_{X_0} \in \tt{\scr P}_p  ($respectively, initial distribution in   $ \tt{\scr P}_p )$,
the SDE $ \eqref{E0} $ has a unique strong $($respectively weak$)$  solution. Moreover,    there exists a constant $c \in (0,\infty)$ such that
\beq\label{*} \|P_t^*\mu \|_{\tt L^{k'}}\le c \| \mu  \|_{\tt L^{p}} t^{-\ff{d(k'-p)}{2pk'}},\ \ \ t\in (0,T],\ \mu\in \tt{\scr P}_p, \ p\le k'\le\infty.\end{equation}
\item[$(2)$] If for some constants $K\in (0,\infty)$  and  $\tau\in [0, \ \infty )$ there holds
\beq\label{P3'} \beg{split}&  |b_t(x,r,\rr)-b_t(x,r',\rr')|\le Kt^\tau \big(|r-r'|+t^{-\ff{d}{2k}} \|\rr-\rr'\|_{\tt L^k}\big),\\
&\qquad\qquad  \ t\in (0,T], \ r,r'\in [0,\infty),\ \rr,\rr'\in \D_1^+,\end{split}\end{equation}
then for any $p'\in [1,p]\cap (\ff d{1+2\tau}\lor \ff{dk }{d+k}, p]$,   there exists a constant $c>0$ such that
\beq\label{R}\beg{split}\|P_t^*\mu- P_t^*\nu\|_{\tt L^k}&\le
 \| \mu  - \nu\|_{\tt L^{p'}} t^{-\ff{d(k-p')}{2p'k}}\exp\Big[c+ct(\|\mu\|_{\tt L^{p'}}\land\|\nu\|_{\tt L^{p'}})^{\ff{2p'}{(1+2\tau)p'-d}} \Big],\\
 &\qquad\qquad  \ \ t\in (0,T],\ \mu,\nu\in \tt{\scr P}_p.\end{split}\end{equation}
In particular, if $d=1$ and   $\tau>0$, then we may take $p'=1$ such that for some   constant $c>0$
\beq\label{R'} \|P_t^*\mu- P_t^*\nu\|_{\tt L^k}\le  c\| \mu  - \nu\|_{\tt L^{1}} t^{-\ff{d(k-1)}{2k}},\ \ \ t\in (0,T],\ \mu,\nu\in \tt{\scr P}_p.\end{equation}
\item[$(3)$] If $\eqref{P3'}$ holds for some $\tau\in [0, \infty)$, then for any $p'\in [1,p]\cap (\ff d{1+2\tau}\lor\frac{dk}{d+k}, k]$ and
\beq\label{PQ}q\in \big[1,\infty\big)\cap \Big(\ff{kd\hat q(p'-1)}{p'[k(\hat q-2)+d\hat q]-kd\hat q},\ \infty\Big), \ \text{where} ~\hat q:= \min_{2\le i\le l} q_i,\end{equation}
there exists a constant $c>0$ such that for any $t\in (0,T]$ and $ \mu,\nu\in \tt{\scr P}_p,$
\beq\label{R1}\beg{split}&\|P_t^*\mu- P_t^*\nu\|_{\tt L^k}\\
&\le \W_q(\mu,\nu)
  t^{-\ff 1 2 -\ff{d}2(\ff 1 {p'}+ \ff{p'-1}{p'q} -\ff 1 k)}\big(\|\mu\|_{\tt L^{p'}}+ \|\nu\|_{\tt L^{p'}}\big)^{\ff {q-1} q}
   \e^{c+ct(\|\mu\|_{\tt L^{p'}}\land\|\nu\|_{\tt L^{p'}})^{\ff{2p'}{(1+2\tau)p'-d}}}.
    \end{split}\end{equation}
If in particular $d=1$ and $\tau >0$, then there exists a constant $c>0$ such that
\beq\label{R1'} \|P_t^*\mu- P_t^*\nu\|_{\tt L^k}\le  c\W_1(\mu,\nu) t^{\ff 1{2k} -1},\ \ \ t\in (0,T],\ \mu,\nu\in \tt{\scr P}_p.\end{equation}\end{enumerate}
\end{thm}

To prove this result,  we will apply the Banach fixed point theorem to  distribution density functions.  We write $f\les g$ for two nonnegative functions if there exists a constant $c\in (0,\infty)$ such that $f\le c g$.

Let $C_w([0,T];\scr P)$ be the class of all weakly continuous maps $\gg: [0,T]\to \scr P$.  For any
\beq\label{GGM} \gg\in {\C}:=\big\{\gg\in C_w([0,T];\scr P):  ~\rr_{\gg(t)}\in  \D_1^+ ,\ t\in (0,T]\big\}, \end{equation} 
 we consider the following SDE with the frozen density parameter $\gg$:
$$ \d X_{s,t}^{\gg,x}= b_t(X_{s,t}^{\gg,x}, \rr_{\gg(t)}(X_{s,t}^{\gg,x}),\rr_{\gg(t)})\d t+\si_t(X_{s,t}^{\gg,x})\d W_t,\ \ s\in [0,T),\ t\in [s,T],\ X_{s,s}^{\gg,x}=x.$$
According to \cite[Proposition 5.1]{HRW25}, the assumption {\bf (A)} implies the well-posedeness of this SDE. Let
$$P_{s,t}^\gg f(x):= \E[f(X_{s,t}^{\gg,x}) ],\ \ \ 0\le s\le t\le T,\ x\in\R^d,\ f\in \B_b(\R^d),$$
where $\B_b(\R^d)$ is the class of bounded measurable functions on $\R^d$. For any $\nu\in \scr P$, let
$P_{s,t}^{\gg*}\nu\in\scr P$ be defined by
$$\big(P_{s,t}^{\gg*}\nu\big)(A):= \int_{\R^d} \big(P_{s,t}^\gg 1_A(x)\big)\nu(\d x)$$ for any measurable set $A\subset\R^d.$
Simply denote
$$P_t^\gg= P_{0,t}^\gg,\ \ \ P_t^{\gg*}= P_{0,t}^{\gg*},\ \ \ t\ge 0.$$
When $s=0$,    the SDE
\beq\label{CP0} \d X_t^\gg= b_t(X_{t}^{\gg}, \rr_{\gg(t)}(X_{t}^{\gg}),\rr_{\gg(t)})\d t+\si_t(X_{t}^{\gg})\d W_t,\ \  \ t\in [0,T],\ X_{0}^{\gg}=X_0\end{equation}
is well-posed.  We define a map $\Phi^\mu=(\Phi^\mu_t)_{t\in [0,T]}: {\C}\to C_w([0,T];\scr P)$ by
\beq\label{CP}  \Phi^\mu_t\gg:= \L_{X_t^\gg}= P_{t}^{\gg*}\mu,\ \ \ t\in [0,T].\end{equation}
 Then
\beq\label{P2'}    \int_{\R^d} P_t^\gg f\d\mu=\E[f(X_t^\gg)]=  \int_{\R^d}(\rr_{\Phi^\mu_t\gg}f)(x) \d x,\ \ \ t\in (0,T],\ f\in\B_b(\R^d),\end{equation}
where the existence of density $\rr_{\Phi^\mu_t\gg}$ follows from the non-degenerate noise, see \cite{BKR}.
Let $\C$ be in \eqref{GGM} and define 
 $${\C}^{p,k}_\mu:=\bigg\{\gg\in {\C}:\ \gg(0)=\mu,\
 \sup_{t\in (0,T]} t^{\ff {d(k-p)}{2pk}}\|\gg(t)\|_{\tt L^k}<\infty\bigg\}.$$    The following lemma implies
 \beq\label{LL1} \Phi^\mu: {\C} \to {\C}^{p,k}_\mu,\ \ \ 1\le p\le k,\  \mu\in \tt{\scr P}_p. \end{equation}
For any $\ll\in (0,\infty)$, \beq\label{LM} d_\ll(\gg,\eta):=\sup_{t\in (0,T]} \Big[\e^{-\ll t} t^{\ff{d(k-p)}{2pk}} \|\gg(t)-\eta(t)\|_{\tt L^k}\Big], \ \ \   \gg,\eta\in {\C}
\end{equation} defines a metric on ${\C}^{p,k}_\mu$ which is complete when $k=1$ but incomplete for $k>1.$

\beg{lem}\label{L1} Assume $\eqref{P1}$ and $\eqref{P2}$ in {\bf (A)}.
Then there exists increasing  $\zeta: [1,\infty)\to  (0,\infty)$
such that
\beq\label{P5} \|\Phi^\mu_t\gg \|_{\tt L^{k_2}}\le \zeta(k_1) \| \mu  \|_{\tt L^{k_1}} t^{-\ff{d(k_2-k_1)}{2k_1k_2}}, \ \ t\in (0,T],\ \gg\in {\C},\ 1\le k_1\le  k_2\le \infty, k_1<\infty.\end{equation}
  \end{lem}

\beg{proof}
   We first introduce some estimates on the regular diffusion semigroup $\bar P_{s,t}$ for the SDE only with drift $b_t^{(1)}$, then estimate $P_{s,t}^\gg$ using
Girsanov's transform and Duhamel's formula, so that \eqref{P5} is proved.

 (a) Consider the reference SDE
\beq\label{E2}\d\bar X_{s,t}^x = b_t^{(1)} (\bar X_{s,t}^{x})\d t+\si_t(\bar X_{s,t}^{x})\d W_t,\ \ s\in [0,T),\ t\in [s,T],\   \bar X_{s,s}^{x}=x.\end{equation}
According to \cite{MPZ}, the conditions on $(a,b^{(1)})$ in {\bf (A)} imply that   the associated semigroup $\bar P_{s,t}$ has heat kernel $\bar p_{s,t}$, i.e.
$$\bar P_{s,t}f(x):= \E[f(\bar X_{s,t}^x)]=\int_{\R^d} \bar p_{s,t}(x,y) f(y)\d y,\ \ \ 0\le s< t\le T,\ x\in\R^d,\ f\in \B_b(\R^d),$$
and for some constants $c,\kk>0$ and diffeomorphism $\psi_{s,t}$ on $\R^d$ with
\beq\label{DD} \|\nn \psi_{s,t}\|+\|\nn\psi_{s,t}^{-1}\| \le c,\end{equation}
the heat kernel $\bar p_{s,t}$ satisfies
\beq\label{DF}  |\nn^i  \bar p_{s,t}(\cdot,y)(x)|\le c (t-s)^{\ff i 2} p_{t-s}^\kk (\psi_{s,t}(x)-y),\ \ \ i=0,1,\
   0\le s<t\le T,\ x,y\in\R^d,\end{equation}
where $\nn^0$ is the identity operator,  $\nn^1$ is the gradient,   and
$$p_t^\kappa(x) = (\pi \kappa t)^{-d/2} \e^{-\frac{|x|^2}{  \kappa t}}, \quad t > 0, \; x \in \mathbb{R}^d
$$ is the  Gaussian heat kernel.
Therefore,  there exist constants $\kk_1,\kk_2\in (0,\infty)$ such that
\beq\label{ES0} \beg{split}&1_{B(x,1)} |\nn^i \bar P_{s,t} (f1_{B(y,1)})|\le  \kk_1 (t-s)^{-\ff i 2} \e^{-\kk_2|x-y|^2} 1_{B(x,1)} P_{t-s}^{2\kk}\big(|f1_{B(y,1)}|\circ \psi_{s,t}\big),\\
& \qquad \ \ i=0,1,\ 0\le s<t\le T,\ x,y\in\R^d,\ f\in \B_b(\R^d),\end{split} \end{equation}
where  $P_{t-s}^{2\kk}f(x):= \int_{\R^d} p_{t-s}^{2\kk}(x-y) f(y)\d y.$
Let $\hat{\Z}^d :=\{d^{-1/2} u: u\in \Z^d\}$.  We find a constant $C(d)\in (1,\infty)$ such that
\beq\label{SD} 1\le \sum_{u\in \hat{\Z}^d} 1_{B(u+z,1)}(x)\le C(d), \ ~x\in \R^d.\end{equation}
Under {\bf (A)}, we have the following  Duhamel's  formula (see \cite[Proposition 5.5]{HRW25})
\beq\label{GF} \beg{split} &P_{s,t}^\gg f= \bar P_{s,t} f + \sum_{i=2}^l\int_s^tP_{s,r}^\gg\big\< b_r^{(i)}(\cdot, \rr_{\gg(r)}(\cdot), \rr_{\gg(r)}),\ \nn \bar P_{r,t}f\big\>\d r,\\
&= \bar P_{s,t} f - \sum_{i=2}^l \int_s^t \bar P_{s,r}\big\<b_r^{(i)}(\cdot, \rr_{\gg(r)}(\cdot), \rr_{\gg(r)}),\ \nn P_{r,t}^\gg f\big\>\d r,\  \ f\in \B_b(\R^d).\end{split}\end{equation}
According to \cite[Lemma   5.3]{HRW25},
\beq\label{5.3}\beg{split}&  \|\nn^i P_{s,t}^{2\kk}\|_{\tt L^{k_1}\to \tt L^{k_2}}+\|\nn^i \bar P_{s,t} \|_{\tt L^{k_1}\to \tt L^{k_2}}\les (t-s)^{-\ff i 2-\ff{d(k_2-k_1)}{2k_1k_2}},\\
&\quad i=0,1,\ 0\le s<t\le T,\ 1\le k_1\le k_2\le \infty.\end{split}\end{equation}

(b) We intend to  find an increasing  function $\tt C: [1,\infty)\to [1,\infty)$ such that for any $\nu\in \scr P$ and $\gg\in {\C}$,
\beq\label{DD0}\beg{split}& \|P_{s,t}^{\gg*}\nu\|_{\tt L^{k_2}}  \le \tt C(k')\|\nu\|_{\tt L^{k_1}} (t-s)^{-\ff{d(k_2-k_1)}{2k_1k_2}},\\
&\quad \ 0\le s<t\le T,\ 1\le k_1\le k_2\le k',\ k'\in [1,\infty),\ \nu\in\scr P,\ \gg\in {\C}.\end{split}\end{equation}
This implies that for any $k'\in  [1,\infty)$,
\eqref{P5} holds for all $1\le k_1\le k_2\le k'$ and $\tt C(k')$ in place of $\zeta(k_1)$.

For fixed $k'\in [1,\infty)$, let   $\theta:= \ff {k'}{k'-1}>1$ depending on $k'$  such that $\ff {k'}{\theta(k'-1)}\ge 1.$ For $1\le k_1\le k_2\le k'$, we have
\beq\label{D1} \ff {k_1}{\theta(k_1-1)}\ge \ff {k_2}{\theta(k_2-1)}\ge \ff {k'}{\theta(k'-1)}\ge  1.\end{equation}
To estimate $ \|P_{s,t}^{\gg*}\nu\|_{\tt L^{k_2}}$,   let
\beq\label{ZK} \D_{z,k_2}:= \Big\{0\le f\in \B_b(\R^d):\ f|_{B(z,1)^c}=0,\ \|f\|_{\ff {k_2}{k_2-1}}\leq 1\Big\},\ \ z\in\R^d.\end{equation}
By \eqref{CP} we have
\beq\label{CP2} \beg{split}& \|P_{s,t}^{\gg*}\nu\|_{\tt L^{k_2}} = \sup_{z\in\R^d} \sup_{f\in \D_{z,k_2}} \nu(P_{s,t}^\gg f)
=  \sup_{z\in\R^d} \sup_{f\in \D_{z,k_2}} \int_{\R^d} \E[f(X_{s,t}^{\gg,x})]\nu(\d x).\end{split}\end{equation}
Let $\bar X_{s,t}^x$ solve the SDE \eqref{E2}, define
\beg{align*} &\xi_t^x:= -(\si_t^*a_t^{-1})(  \bar X_{s,t}^{x} ) \bigg(\sum_{i=2}^l b_t^{(i)}(  \bar X_{s,t}^{x} , \rr_{\gg(t)}(  \bar X_{s,t}^{x} ), \rr_{\gg(t)})\bigg),\\
& \tt W_t:=W_t- \int_s^t \xi_r^x\d r,\ \ \ t\in [s,T]. \end{align*}
By \eqref{P1}, \eqref{P2}  and Khasminski's estimate, see \cite[(1.2.7) and Theorem 1.2.4]{RW25} for $B_s=0$, for any $\ll\in (0,\infty)$ we find a constant $c(\ll)\in (0,\infty)$ such that
$$\E \big[\e^{\ll \int_0^T |\xi_s^x|^2\d s}\big] \le c(\ll),\ \ x\in \R^d.$$
So, by Girsanov's theorem, $\tt W_t$ is an $m$-dimensional Brownian motion under the probability measure $\Q:=  R_{s,x} \P$, where
$$R_{s,x}:= \e^{\int_s^T \<\xi_r^x,\d W_r\>-\ff 1 2 \int_s^T |\xi_r^x|^2\d r}$$
satisfies
$$h(n):= \sup_{x\in\R^d, s\in [0,T)} \big(\E[R_{s,x}^{n}]\big)^{\ff 1 n} <\infty,\ \ \ n\in [1,\infty),$$
which is increasing in $n$. 
Combining this with the weak uniqueness of the SDE for $\bar X_{s,t}^x$,   \eqref{CP2}, and recalling $\theta=\ff{k'}{k'-1}$,  we
obtain
\beq\label{CP3} \beg{split}&\|P_{s,t}^{\gg*}\nu\|_{\tt L^{k_2}} = \sup_{z\in\R^d} \sup_{f\in \D_{z,k_2}} \int_{\R^d} \E[ R_{s,x}  f(\bar X_{s,t}^x)]\nu(\d x)\\
&\le  \sup_{z\in\R^d} \sup_{f\in \D_{z,k_2}} \int_{\R^d} {\big(\E[R_{s,x}^{k'}]\big)^{\ff {1} {k'} }}  \big(\E[f(\bar X_{s,t}^x)^\theta]\big)^{\ff 1 \theta}\nu(\d x)\\
&\le h(k') \sup_{z\in\R^d} \sup_{f\in \D_{z,k_2}} \int_{\R^d} \big(\E[ f(\bar X_{s,t}^x)^{\theta}]\big)^{\ff 1 \theta}\nu(\d x)\\
&= h(k') \sup_{z\in\R^d} \sup_{f\in \D_{z,k_2}}\nu\Big(\big(\bar P_{s,t } f^\theta \big)^{\ff 1 \theta}\Big).\end{split}
\end{equation}
  By \eqref{ES0}, \eqref{SD}, \eqref{5.3} and \eqref{D1} for $i=0$, we find constants $c_1,c_2 \in (0,\infty)$ such  that for any $f\in \D_{z,k_2}$,
 {\beg{align*} &\nu\Big(\big(\bar P_{s,t} f^\theta \big)^{\ff 1 \theta}\Big)
\le  \sum_{u\in \hat{\Z}^d}  \nu\Big(1_{B(u+z,1)}\big(\bar P_{t-s}    f^\theta   \big)^{\ff 1 \theta}\Big) \le  \sum_{u\in \hat{\Z}^d}  \|\nu\|_{\tt L^{k_1}} \big\| \bar P_{t-s}  f^\theta \big\|_{\tt L^{\ff{k_1}{\theta(k_1-1)}}}^{\ff 1 \theta}\\
&\le c_1 \|\nu\|_{\tt L^{k_1}} \big\| \bar P_{t-s} \big\|_{{\tt{L}}^{\ff{k_2}{\theta(k_2-1)}}\to  {\tt{L}}^{\ff{k_1}{\theta( k_1-1)}}}^{\ff 1 \theta} { \|f^\theta   \big\|_{\tt L^{\ff{k_2}{\theta(k_2-1)}}}^{\ff 1 \theta}}\\
&\le c_2\|\nu\|_{\tt L^{k_1}} (t-s)^{-\ff{d(k_2-k_1)}{2k_1k_2}} { \|f   \big\|_{\tt L^{\ff{k_2}{ k_2-1}}} },\ \ 1\le k_1\le k_2\le k',\ 0\le s<t\le T.\end{align*}
}
Combining this with \eqref{CP3} we prove \eqref{DD0} for some $\tt C(k')\in (0,\infty)$ increasing in $k'$.

(c)  Since $(p_i,q_i)\in \scr K,$ we fixed  $k'\in (d,\infty)$ such that $k'>\hat p:= \sup_{2\le i\le l}\ff{p_i}{p_i-1}$ and
$$\vv:= \max_{2\le i\le l}\Big(\ff 1 2 +\ff d {2k'}+\ff d {2p_i}\Big)\ff{q_i}{q_i-1}<1.$$
By \eqref{DD0} , it remains to find  an increasing function   $\zeta: [1,\infty)\to  (0,\infty)$ such that
\eqref{P5} holds for all $k_2\ge  k'$ and $k_1\in [1,k_2]\cap [1,\infty)$.

 We first consider the case that
 $k_1\ge  k'.$   In this case, we have $k_1\ge \ff{p_i}{p_i-1}$ since $k'\ge \hat p,$ and
 $$\Big[\ff 1 2 +\ff d 2\Big(\ff 1 {k_1}+\ff 1 {p_i}-\ff 1 {k_2}\Big)\Big]\ff{q_i}{q_i-1}\le \vv<1.$$
So,  we find  constants $c_1,c_2\in (0,\infty)$ such that \eqref{5.3}   yields
\beg{align*}   &\sum_{u\in \hat{\Z}^d}  \bigg(\int_0^t \|1_{B(u+z,1)} \nn\bar P_{t-s}f\|_{\ff{p_ik_1}{p_i(k_1-1)-k_1}}^{\ff {q_i}{q_i-1}}  \d s\bigg)^{\ff{q_i-1}{q_i}}\\
&\le  c_1   \bigg(\int_0^t\Big[ (t-s)^{-\ff 1 2-\ff d 2 (\ff 1{k_1}+ \ff 1 {p_i} -\ff 1 {k_2})} \Big]^{\ff {q_i}{q_i-1}}\d s\bigg)^{\ff{q_i-1}{q_i}}\\
&\le c_2 t^{-\ff{d(k_2-k_1)}{2k_1k_2}}, \ \ \  t\in (0,T],\ z\in\R^d,\ f\in \D_{z,k_2}.\end{align*}
Combining this with   H\"older's inequality, we find a constant $c_3\in (0,\infty)$ such that
\beg{align*} &\sum_{2\le i\le l}\sum_{u\in \hat{\Z}^d} \int_0^t \big\|1_{B(u+z,1)} f_s^{(i)} {\nn\bar P}_{t-s} f\big\|_{\ff{k_1}{k_1-1}} \d s\\
&\le \sum_{2\le i\le l}  \int_0^t\sum_{z\in \hat{\Z}^d}\big\|1_{B(u+z,1)} f_s^{(i)}\big\|_{p_i}  \big\|1_{B(u+z,1)}|\nn\bar P_{t-s} f|\big\|_{\ff{p_ik_1}{p_i(k_1-1)-k_1}}\d s\\
&\le \sum_{2\le i\le l} \|f^{(i)}\|_{\tt L_{q_i}^{p_i}}  \sum_{u\in \hat{\Z}^d} \bigg(\int_0^t \big\|1_{B(u+z,1)}|\nn\bar P_{t-s} f|\big\|_{\ff{p_ik_1}{p_i(k_1-1)-k_1}}^{\ff{q_i}{q_i-1}}\d s\bigg)^{\ff{q_i-1}{q_i}}\\
&\le c_3 t^{-\ff{d(k_2-k_1)}{2k_1k_2}}, \ \ \  t\in (0,T],\ z\in\R^d,\ f\in \D_{z,k_2}. \end{align*}
 Noting that \eqref{DD0} implies  $\|\Phi_s^\mu\gg\|_{\tt L^{k_1}}\le \tt C(k_1) \|\mu\|_{\tt L^{k_1}},$
 combining this with \eqref{P2}  and \eqref{GF}, we find a constant $c_4\in (0,\infty)$ increasing in  $k_1$ such that
 for any $k_2\ge k_1\ge k'$ and $k_1<\infty$,
\beg{align*} &\mu(P_{t}^\gg f) \le \mu(\bar P_{t}f) + \sum_{i=2}^l\sum_{u\in \hat{\Z}^d}   \int_0^t \|\Phi_{s}^\mu\gg\|_{\tt L^{k_1}}\|1_{B(u+z,1)}f_s^{(i)}{\nn \bar P}_{t-s}  f\|_{\ff {k_1}{k_1-1}}\d s\\
&\le c_4  t^{-\ff{d(k_2-k_1)}{2k_1k_2}}{ \|\mu\|_{\tt L^{k_1}}} ,\ \ k'\le k_1\le k_2,\ t\in (0,T],\ f\in  \D_{z,k_2},\ z\in\R^d.\end{align*}
So, \eqref{P5} holds for some increasing $\zeta: [1,\infty)\to  (0,\infty)$ and all $k_2\ge k_1 \ge k', k_1<\infty.$

Next, let $k_1\le k'\le k_2$.  By   the semigroup property
$$\Phi_t^\mu\gg = P_t^{\gg*}\mu= P_{\ff t 2,t}^{\gg*} P_{\ff t 2}^{\gg*} \mu= P_{\ff t 2,t}^{\gg*}\Phi_{\ff t 2}^\mu\gg,$$
   \eqref{P5}  for $k_2\ge k_1=k'$ which has been just proved, and  \eqref{DD0} for $k_2=k'\ge k_1$ which has been proved in step (b),  we find a constant  $c_5\in (0,\infty)$ increasing in 
    $k_1$ such that
\beg{align*} &\|\Phi_t^\mu\gg\|_{\tt L^{k_2}}= \sup_{z\in \R^d}\sup_{f\in \D_{z,k_2}} |(\Phi_{t/2}^\mu\gg)(P_{t/2,t}^\gg f)|\\
&\le \tt C(k')\|\Phi_{\ff t 2}^\mu\gg\|_{\tt L^{k'}} (t/2)^{-\ff{d(k_2-k')}{2k'k_2}}\le c_5 \|\mu\|_{\tt L^{k_1}} t^{-\ff{d(k_2-k_1)}{2k_1k_2}},\ \ 1\le k_1\le  k'\le k_2.\end{align*}
  Then the proof is finished.
 \end{proof}

By  \eqref{LL1}, any solution of \eqref{E0} with initial distribution $\mu$ satisfies $(\L_{X_t})_{t\in [0,T]}\in \C^{p,k}_\mu.$ Combining this with the well-posedness of \eqref{CP0},
to show that \eqref{E0} has a unique (weak/strong) solution for $\mu\in \tt{\scr P}_p$, we only need to prove  that $\Phi^\mu$ has a   unique fixed point in ${\C}^{p,k}_\mu$.
To this end, we show that when $p\in (d,k]\cap [1,\infty)$, $\Phi^\mu:  {\C}_\mu^{p,k}\to {\C}_\mu^{p,k} $ is contractive under the metric $d_\ll$ defined in \eqref{LM} for  some $\ll\in (0,\infty)$.

\beg{lem}\label{L2} Assume {\bf (A)} and let $p\in (d,k]\cap [1,\infty)$. Then there exists a constant $\ll\in (0,\infty)$ such that
$$ d_\ll(\Phi^\mu\gg, \Phi^\mu\eta)\le \ff 1 2 d_\ll(\gg,\eta),\ \ \ \gg,\eta\in {\C}.$$
\end{lem}

\beg{proof} To estimate $\|\Phi^\mu_t\gg-\Phi^\mu_t\eta\|_{\tt L^k}$ for $\gg,\eta\in {\C},$
we take  $f\in \D_{z,k}$ in \eqref{ZK}  for $z\in\R^d$.    By  \eqref{GF},  {\bf (A)} and \eqref{ES0}, we obtain
\beq\label{C1}  \beg{split} & |(\Phi^\mu_t\gg)(f)- (\Phi^\mu_t\eta)(f)|= |\mu(P_t^\gg f-P_t^\eta f)| \\
 &=  \sum_{i=2}^l \int_0^t \Big\| \big\<\rr_{\Phi^\mu_s\gg}  b_s^{(i)} (\cdot, \rr_{\gg(s)}(\cdot),\rr_{\gg(s)})- \rr_{\Phi^\mu_s\eta}
   b_s^{(i)} (\cdot, \rr_{\eta(s)}(\cdot,\rr_{\eta(s)}), \nn \bar P_{s,t} f\big\>\Big\|_1\d s\\
 &\les   \sum_{i=2}^l  \int_0^t (t-s)^{-\ff 1 2} \Big(\Big\| \big(\rr_{\Phi^\mu_s\gg} -\rr_{\Phi^\mu_s\eta}\big)   b_s^{(i)} (\cdot, \rr_{\eta(s)}(\cdot),\rr_{\eta(s)})   P_{s,t}^{2\kk} (f\circ \psi_{s,t})     \Big\|_{1}\\
 &\qquad\qquad +\Big\| \rr_{\Phi^\mu_s\gg}  \big[b_s^{(i)} (\cdot, \rr_{\gg(s)}(\cdot),\rr_{\gg(s)})-    b_s^{(i)}  (\cdot, \rr_{\eta(s)}(\cdot),\rr_{\eta(s)})\big]
    P_{s,t}^{2\kk} (f\circ \psi_{s,t}) \Big\|_{1}\Big) \d s\\
 &\les  \sum_{i=1}^l \int_0^t (t-s)^{-\ff 1 2} B_i(s)\d s,\end{split}\end{equation}
 where
 \beg{align*} &  B_1(s):=    \Big\|   \rr_{\Phi^\mu_s\gg} \big(|\rr_{\gg(s)}- \rr_{\eta(s)}| +s^{-\ff d{2k}}\|\gg(s)-\eta(s)\|_{\tt L^k}\big)  P_{s,t}^{2\kk} (f\circ \psi_{s,t}) \Big\|_{1},\\
 &B_i(s):=  \Big\|  \big(\rr_{\Phi^\mu_s\gg}-  \rr_{\Phi^\mu_s\eta} \big) f_s^{(i)}  P_{s,t}^{2\kk} (f\circ \psi_{s,t}) \Big\|_{1},\ \ \ 2\le i\le l.\end{align*}
 By $f\in \D_{z,k}$, $\eqref{DD}$,   the symmetry of $P_{s,t}^{2\kk}$ in $L^2(\d x)$,   \eqref{5.3}, \eqref{P5}, and H\"older's inequality, we have
 \beg{align*}& \int_0^t (t-s)^{-\ff 1 2} B_1(s) \d s  \\
 &\les  \int_0^t (t-s)^{-\ff 1 2} \Big\|  P_{t-s}^{2\kk}\big[\rr_{\Phi^\mu_s\gg} \big(|\rr_{\gg(s)}- \rr_{\eta(s)}| +s^{-\ff d{2k}}\|\gg(s)-\eta(s)\|_{\tt L^k}\big) \Big\|_{\tt L^k}\d s \\
 &\les  \int_0^t (t-s)^{-\ff 1 2} \|P_{t-s}^{2\kk}\|_{\tt L^{p} \to \tt L^k} \Big(  \| \Phi^\mu_s\gg \|_{\tt L^{\ff {kp}{k-p}}}+s^{-\ff d{2k}} \|\Phi^\mu_s\gg \|_{\tt L^{p}}\Big)
 \| \gg(s)-  \eta(s)\|_{\tt L^k}\d s \\
 &\les  \| \mu \|_{\tt L^p} \int_0^t  (t-s)^{-\ff 1 2 -\ff{d(k-p)}{2kp}} s^{-\ff d {2k}}   \| \gg(s)-  \eta(s)\|_{\tt L^k}\d s,\ \ t\in [0,T].   \end{align*}
 Similarly,  for $2\le i\le l,$ by \eqref{ES0},    H\"older's inequality and \eqref{SD},    we obtain
\beg{align*} & \int_0^t (t-s)^{-\ff 1 2}B_i(s)\d s\\
 &\les \sum_{u\in\hat\Z^d} \e^{-\kk_2|u|^2}  \int_0^t(t-s)^{-\ff 1 2} \Big\|  P_{t-s}^{2\kk} \big[1_{B(u+ z,1) }   |\rr_{\Phi^\mu_s\gg}-  \rr_{\Phi^\mu_s\eta}|
 f_s^{(i)} ]  \Big\|_{k}\d s\\
&\les \sum_{u\in\hat\Z^d} \e^{-\kk_2|u|^2} \int_0^t(t-s)^{-\ff 1 2} \|P_{t-s}^{2\kk}\|_{\tt L^{\ff{p_ik}{p_i+k}} \to \tt L^k} \|1_{B(u+
 z,1)}    f^{(i)}_s\|_{p_i} \| \Phi^\mu_s\gg -   \Phi^\mu_s\eta \|_{\tt L^k} \d s \\
&\les  \|f^{(i)}\|_{\tt L^{p_i}_{q_i}(0,t)} \bigg(\int_0^t \Big[(t-s)^{-\ff 1 2 -\ff d{2p_i}}    \| \Phi^\mu_s\gg -   \Phi^\mu_s\eta \|_{\tt L^k}\Big]^{\ff{q_i}{q_i-1}}\d s\bigg)^{\ff{q_i-1}{q_i}},\ \ t\in [0,T],  \ 2\le i\le l.\end{align*}
Combining these with \eqref{C1} and $f^{(i)}\in \tt L_{q_i}^{p_i}$, we derive
 \beg{align*} &\|\Phi^\mu_t\gg - \Phi^\mu_t\eta\|_{\tt L^k}= \sup_{z\in \R^d}\sup_{f\in {\D}_{z,k}}  |(\Phi^\mu_t\gg)(f)- (\Phi^\mu_t\eta)(f)|\\
&\les   \|\mu\|_{\tt L^p}  \int_0^t (t-s)^{-\ff 1 2-\ff{d(k-p)}{2kp}} s^{-\ff d {2k}  }   \| \gg(s)-  \eta(s)\|_{\tt L^k}\d s \\
&\quad + \sum_{i=2}^l { \|f^{(i)}\|_{\tt L^{p_i}_{q_i}(0,t)}} \bigg(\int_0^t \Big[(t-s)^{-\ff 1 2 -\ff d{2p_i}}    \| \Phi^\mu_s\gg -   \Phi^\mu_s\eta \|_{\tt L^k}\Big]^{\ff{q_i}{q_i-1}}\d s\bigg)^{\ff{q_i-1}{q_i}},\ \ t\in [0,T].\end{align*}
So, we find a constant $C_1>0$ such that
\beq\label{C*} \beg{split} &d_\ll(\Phi^\mu\gg,\Phi^\mu\eta)\\
&\le C_1  d_\ll(\gg,\eta)   \sup_{t\in (0,T] } t^{\ff{d(k-p)}{2pk}}
\int_0^t (t-s)^{-\ff 1 2-\ff {d(k-p)}{2pk}}s^{-\ff d{2p}} \e^{-\ll(t-s)}  \d s\\
&\quad + C_1 d_\ll(\Phi^\mu\gg,\Phi^\mu\eta) \sup_{t\in (0,T] } t^{\ff{d(k-p)}{2pk}}
\sum_{i=2}^l \bigg(\int_0^t (t-s)^{-\vv_i} s^{-\dd_i}\e^{-   \ll (t-s)} \d s\bigg)^{\ff{q_i-1}{q_i}},\end{split}\end{equation}
where, by  $(p_i,q_i)\in \scr K$ and $k\ge p>d$,
\beq\label{AE} \beg{split}& \vv_i:= \Big(\ff 1 2+ \ff d{2p_i}\Big)\ff{q_i}{q_i-1}<1,\ \  \ \
\dd_i:= \ff{d(k-p)q_i}{2pk(q_i-1)} <1,\\
&\ff 1 2 +\ff{d(k-p)}{2pk} <1,\ \ \ \  \ff 1 2- \ff d{2p}>0,\\
&\ff{d(k-p)}{2pk}+ \ff{q_i-1}{q_i} \big(1-\vv_i-\dd_i\big)= \ff 1 2 -\ff 1 {q_i} -\ff d{2p_i}>0. \end{split}  \end{equation}
So,
\beq\label{AD}  \theta:= \min_{2\le i\le l}\Big\{ \ff{d(k-p)}{2pk} - \ff{q_i-1}{q_i} \big(\vv_i{+}\dd_i-1\big),\ \ff 1 2- \ff{d}{2p}\Big\}\in \Big(0,\ff 1 2\Big).\end{equation}
By the FKG inequality,   for any $\vv,\dd\in (0,1)$ with   $\theta\in (0,1-\vv)$,  we find a constant $c>0$ such that
\beq\label{FKG} \beg{split} &\int_0^t (t-s)^{-\vv}s^{-\dd} \e^{-\ll (t-s)}\d s\le \bigg(\ff 1t \int_0^t s^{-\dd}\d s\bigg)\int_0^t   (t-s)^{-\vv}\e^{-\ll (t-s)}\d s\\
&\le \ff {t^{-\dd}} {1-\dd}  \int_0^t t^{1-\vv-\theta}  (t-s)^{\theta-1}\e^{-\ll (t-s)}\d s \le c    t^{1-\vv-\dd-\theta} \ll^{-\theta},\ \ \ t\in (0,T],\ \ll>0.\end{split}\end{equation}
 Combining this with \eqref{C*} and \eqref{AD},
  we find constants $C_1,C_2>0$  possibly depending on $\|\mu\|_{\tt L^p}$,     such that
\beg{align*} d_\ll(\Phi^\mu\gg,\Phi^\mu\eta) \le C_2 d_\ll( \gg ,  \eta  ) \ll^{-\theta},\ \ \ll>0.\end{align*}
 So, the desired estimate holds for  $\ll\ge (2C_2)^{1/\theta}$.
\end{proof}

\beg{proof}[Proof of Theorem \ref{T1}] Assume {\bf (A)} and  let  $p\in (d,k]\cap [1,\infty).$

(1) Let $\mu\in \tt{\scr P}_p$.  For any solution $X_t$ of \eqref{E0} with initial distribution $\mu$,
 $\gg(t):= \L_{X_t}$ for $t\in [0,T]$ satisfies  $\Phi^\mu\gg=\gg\in {\C},$ so that \eqref{*} follows from   Lemma \ref{L1}.
 Hence, the first assertion follows if the map $\Phi^\mu$ has a unique fixed point in   ${\C}_\mu^{p,k}$.  In this case, the unique solution
  is $X_t^\gg$ for the unique fixed point $\gg$ of $\Phi^\mu$.

 By Lemma \ref{L2}, $\Phi^\mu$ has at most one fixed point in ${\C}_\mu^{p,k}$. However, since ${\C}_\mu^{p,k}$ is not complete under $d_\ll$, Lemma \ref{L2} is not  enough to imply the existence of fixed point.
 To overcome this problem,  we extend ${\C}_\mu^{p,k}$ to a complete space ${\C}_{\mu,sb}^{p,k}$.
  More precisely, let
 $\scr P_{sb}$ be the space of sub-probability measures on $\R^d$, denote again $\rr_\mu(x):=\ff{\mu(\d x)}{\d x}$ for absolutely continuous $\mu\in \scr P_{sb}$, and let $C_w([0,T]; \scr P_{sb})$ be the space of weakly continuous maps from $[0,T]$ to $\scr P_{sb}$. Then
  $${\C}_{\mu,sb}^{p,k}:= \bigg\{\gg\in C_w([0,T]; \scr P_{sb}):\ \gg(0)=\mu, \sup_{t\in (0,T]}   t^{\ff {d(k-p)}{2pk}}
  \|\rr_{\gg(t)}\|_{\tt L^k}<\infty\bigg\}$$
  is complete under the metric $d_\ll$ for any $\ll\in (0,\infty)$, which is defined in \eqref{LM} for $\xi,\eta\in {\C}_{\mu,sb}^{p,k}.$ Now, for fixed $\gg^{(1)}\in {\C}_\mu^{p,k}$,
  let
  $$\gg^{(n+1)}:= \Phi^\mu \gg^{(n)},\ \ \ n\ge 1.$$
  By Lemma \ref{L2}, there exists  $\ll\in (0,\infty)$ such that
  $$  d_\ll(\gg^{(n+2)},\gg^{(n+1)})\le 2^{-n}  d_\ll(\gg_2,\gg_1),\ \ \ n\ge 2,$$ so that
 $\{\gg^{(n)}\}_{n\ge 1}$ is a Cauchy sequence under $ d_\ll$. Then there exists a unique  $\gg^{(\infty)}\in {\C}_{\mu,sb}^{p,k}$ such that
 \beq\label{99} \lim_{n\to\infty}  d_\ll({\gg^{(n)}},\gg^{(\infty)})=0.\end{equation}
  It remains to show that $\gg^{(\infty)}\in   {\C}_{\mu}^{p,k}$, which together with Lemma \ref{L2} and \eqref{99} implies that $\gg^{(\infty)}$ is a  fixed point of $\Phi^\mu$ in   $\C_{\mu}^{p,k}.$ 
 By  \cite[Theorem 1.3.1]{RW25},    {\bf (A)} implies
  $$\E\bigg[\sup_{t\in (0,T]}|X_t^{\gg }|^2\bigg|\F_0\bigg]\le c(1+|X_0|^2),\ \ {\gg }\in {\C}    $$ for some constant $c\in (0,\infty)$. Then
  \beg{align*}& \lim_{N\to\infty} \sup_{{\gg }\in\C,t\in (0,T]} \P(|X_t^{\gg }|^2\ge N) =\lim_{N\to\infty} \sup_{{\gg }\in\C,t\in (0,T]}\E\big[\E(1_{\{|X_t^{\gg }|^2\ge N\}}|\F_0)\big]\\
  &\le \lim_{N\to\infty} \sup_{{\gg }\in\C,t\in (0,T]}\E\bigg[ 1\land \E\Big(\ff{|X_t^{\gg }|^2 }N \Big|\F_0\Big)\bigg]\le \lim_{N\to\infty} \E\bigg[ 1\land \ff{c(1+|X_0|^2)}N\bigg]\\
  &\le \lim_{N\to\infty}  \Big[\vv+ \P\big( c(1+|X_0|^2) >\vv N\big)\Big] =\vv,\ \ \vv>0,\ t\in (0,T].\end{align*}
  Consequently,    for any $t\in (0,T]$,
  $\{\gg^{(n)}(t)\}_{n\ge 1}$ is tight, which together with \eqref{99} implies $\gg^{(\infty)}(t)(\R^d)=1$.   Combining this with  $\gg^{(\infty)}\in {\C}_{\mu,sb}^{p,k}$ we prove  $\gg^{(\infty)}\in   {\C}_{\mu}^{p,k}.$

 (2)  Let \eqref{P3'} hold for some $\tau\in [0, \infty ),$  let $p'\in  [1,p]\cap (\ff d{1+2\tau}\lor \ff{dk }{d + k }, p]$ and $\mu,\nu\in \tt{\scr P}_p.$   We intend to prove \eqref{R}.
 Without loss of generality, we assume that $\|\mu\|_{\tt L^{p'}}\le \|\nu\|_{\tt L^{p'}}$. Denote
 $$\gg(t)= P_t^*\mu, \ \ \ \eta(t)=   P_t^*\nu,\ \ \ t\in [0,T].$$
 Since $\gg$ and $\eta$ are fixed points of   $\Phi^\mu$   on ${\C}_{\mu}^{p,k}$ and $\Phi^\nu$ on ${\C}_\nu^{p,k}$ respectively, by \eqref{P2} and \eqref{GF}, we have
\beq\label{O0}  \beg{split} & \big|(P_t^*\mu)(f)- (P_t^*\nu)(f)\big| \leq \big|(\mu-\nu)(\bar P_t f)\big|\\
 & +\sum_{2\le i\le l} \int_0^t  \Big\|\Big\<\rr_{\gg(s)} b_s^{(i)}(\cdot,\rr_{\gg(s)}(\cdot), \rr_{\gg(s)})- \rr_{\eta(s)}b_s^{(i)}(\cdot,\rr_{\eta(s)}(\cdot), \rr_{\eta(s)}),
 \ \nn \bar P_{s,t}f\Big\>\Big\|_1 \d s.
\end{split}\end{equation}
By \eqref{DD}, \eqref{DF}, the symmetry of $P_t^{2\kk}$  and \eqref{5.3},
\beq\label{89}\beg{split}& \sup_{z\in \R^d}\sup_{f\in \D_{z,k}} |(\mu-\nu)(\bar P_t f)| \les  \sup_{z\in \R^d}\sup_{f\in \D_{z,k}}  \big\|   |\rr_{\mu}-\rr_\nu|  P_t^{2\kk} (f\circ \psi_{0,t})  \big\|_{1}\\
&= \sup_{z\in \R^d}\sup_{f\in \D_{z,k}}  \big\|   (f\circ \psi_{0,t})   P_t^{2\kk}   |\rr_{\mu}-\rr_\nu| \big\|_{1} \les \big\|P_t^{2\kk}   |\rr_{\mu}-\rr_\nu| \big\|_{ \tt L^k}\\
&\les t^{-\ff {d(k-p')}{2p'k}}\|\mu-\nu\|_{\tt L^{p'}},     \ \ \ t\in (0,T],\ \mu,\nu\in \tt{\scr P}_p,\ p'\in [1,p].\end{split} \end{equation}
Moreover, by using \eqref{P3'} and $p'$ in place of \eqref{P3} and $p$ respectively,  the same argument leading to \eqref{C*} yields
\beg{align*} & \sum_{2\le i\le l} \int_0^t   \Big\|\Big\<\rr_{\gg(s)} b_s^{(i)}(\cdot,\rr_{\gg(s)}(\cdot), \rr_{\gg(s)})- \rr_{\eta(s)}b_s^{(i)}(\cdot,\rr_{\eta(s)}(\cdot), \rr_{\eta(s)}),
\ \nn \bar P_{s,t}f\Big\> \Big\|_1  \d s\\
&\les  \|\mu\|_{\tt L^{p'}}  \int_0^t (t-s)^{-\ff 1 2  -\ff{d(k-p')}{2p'k}} s^{\tau -\ff d{2k}}  \| \gg(s) -   \eta(s) \|_{\tt L^k}\d s\\
& \quad  +  \sum_{i=2}^l { \|f^{(i)}\|_{\tt L^{p_i}_{q_i}(0,t)}} \bigg(\int_0^t \big[ (t-s)^{-\ff 1 2 -\ff d  {2p_i}} \| \gg(s) -   \eta(s) \|_{\tt L^k}\big]^{\ff{q_i}{q_i-1}}d s\bigg)^{\ff{q_i-1}{q_i}}, \\
&\qquad\qquad  \  t\in [0,T],\ \mu,\nu\in \tt{\scr P}_p. \end{align*}
Combining this with   \eqref{O0} and \eqref{89},  we find a constant $C_1\in (0,\infty)$ such that
\beq\label{90} \beg{split} &\| \gg(t) - \eta(t)\|_{\tt L^k}\le C_1  t^{-\ff {d (k-p')}{2p'k}}  \|\rr_\mu-\rr_\nu\|_{\tt L^{p'}} \\
&+  C_1 \|\mu\|_{\tt L^{p'}}  \int_0^t (t-s)^{-\ff 1 2  -\ff{d(k-p')}{2p'k}} s^{\tau -\ff d{2k}}  \| \gg(s) -   \eta(s) \|_{\tt L^k}\d s\\
&  + C_1\sum_{i=2}^l { \|f^{(i)}\|_{\tt L^{p_i}_{q_i}(0,t)}} \bigg(\int_0^t \big[ (t-s)^{-\ff 1 2 -\ff d  {2p_i}} \| \gg(s) -   \eta(s) \|_{\tt L^k}\big]^{\ff{q_i}{q_i-1}}d s\bigg)^{\ff{q_i-1}{q_i}},  \ \ t\in [0,T]. \end{split}\end{equation}
By  $(p_i,q_i)\in \scr K$ and  $ p'>\ff d{1+2\tau}\lor \frac{dk}{d+k}$, we have
\beg{align*}&\vv_i:= \Big(\ff 1 2+ \ff d{2p_i}\Big)\ff{q_i}{q_i-1}<1,\ \ \ \dd_i':= \ff{d(k-p')q_i}{2p'k(q_i-1)}<1,\\
&  \ff 1 2 +\ff{d(k-p')}{2p'k}<1,\ \ \   \ \ \ff 1 2 +\ff d{2p'} -\tau<1.\end{align*}
Moreover, since   $\gg=\Phi^\mu\gg,  \eta=\Phi^\nu \eta$,     \eqref{*} implies
$$I_\ll :=\sup_{t\in (0,T]} \|\rr_{\gg(t)} - \rr_{\eta(t)}\|_{\tt L^k} t^{\ff {d(k-p')}{2p'k}}   \e^{-\ll t}<\infty,\ \ \ \ll\ge 0.$$
So, for
$$ \theta_1:= \min_{2\le i\le l}\Big\{\ff 1 2 -\ff 1 {q_i}-\ff d{2p_i}\Big\}, \ \ \ \theta_2:=\ff{(1+2\tau)p'-d}{2p'},$$
 we find a constant $C_2>0$ such that \eqref{90} implies
\beg{align*}I_\ll &\le C_1\|\mu-\nu\|_{\tt L^{p'}} + I_\ll C_1    \|\mu\|_{\tt L^{p'}}
 \sup_{t\in (0,T]}t^{\ff {d(k-p')}{2p'k}}  \int_0^t
 (t-s)^{-\ff 1 2  -\ff{d(k-p')}{2p'k}} s^{\tau -\ff{d}{2p'}} \e^{-\ll(t-s)}\d s\\
&\quad +I_\ll C_1  \sup_{t\in (0,T]}t^{\ff {d(k-p')}{2p'k}} \sum_{i=2}^l
\left(\int_0^t (t-s)^{-\vv_i } s^{-\dd_i'}\e^{-\ll (t-s)}\d s \right)^{\ff{q_i-1}{q_i}}\\
&\le C_1\|\mu-\nu\|_{\tt L^{p'}} + I_\ll C_2 \big(\|\mu\|_{\tt L^{p'}}\ll^{-\theta_2} +\ll^{-\theta_1}\big).\end{align*}
As we have assumed $\|\mu\|_{\tt L^{p'}}\le \|\nu\|_{\tt L^{p'}},$ we prove \eqref{R} for some constant $C>0$ by taking
 $$\ll= \big(4 C_2\|\mu\|_{\tt L^{p'}}\big)^{\theta_2^{-1}} + (4C_2)^{\theta_1^{-1}}.$$
 Note that $\|\mu\|_{\tt L^1}\le 1$. When $d=1$ and $\tau>0$, we may take $p'=1$, so that \eqref{R} implies \eqref{R'}.

 (3)  Comparing with the proof \eqref{R}, the key different in the proof of \eqref{R1} is to derive an alternative   estimate to \eqref{89} by using $\W_q(\mu,\nu)$ in place
 of $\|\mu-\nu\|_{\tt L^{p'}}$. To this end, we recall the maximal functional for a nonnegative measurable function $g$:
 $$\scr M g(x):= \sup_{r\in (0,1)} \ff 1 {|B(x,r)|}\int_{B(x,r)} g(y)\d y,\ \ x\in\R^d.$$
 Then there exists a constant $c_1>0$ such that for any $0\le g\in C_b^1(\R^d),$
 \beq\label{01}\beg{split} & |  \bar P_tg(x)-  \bar P_tg(y)|\le c_1 |x-y|\big(\scr M   |\nn\bar  P_tg|(x)+\scr M | \nn \bar P_tg|(y) + \|\bar P_tg\|_{\infty}  \big), \\
 &\|\scr M  |\nn \bar  P_tg|\|_{\tt L^n}\le c_1 \| \nn \bar  P_tg\|_{\tt L^n},\ \ \ n\in [1,\infty],\end{split}\end{equation} see \cite[Lemma 2.1]{XXZZ}.
 Let $q$ be in \eqref{PQ}, so that
  \beq\label{dd'}\dd':= \ff 1 2 + \ff d 2\Big(\ff 1 {p'}+\ff{p'-1}{qp'}-\ff 1 k\Big)<\ff{\hat q-1}{\hat q}=\min_{2\le i\le l}\ff{q_i-1}{q_i}.\end{equation}
  Then $k\ge d  \lor p' $  implies  $\dd' \ge \ff 1 2\ge \ff{d}{2k}.$
 By this and \eqref{5.3} for  $i=0, k_2=\infty$ and $k_1=k$,  we derive
 \beq\label{04} \|\bar P_tf\|_\infty\les t^{-\ff d{2k}} \les t^{-\dd'},\ \ z\in \R^d,\ f\in \D_{z,k},\ t\in (0,T].\end{equation}
Combining this with \eqref{01} and H\"older's inequality, we find a constant $  c_2 \in (0,\infty)$ such that   for any $z\in \R^d$ and $f\in \D_{z,k},$
 \beq\label{02}\beg{split} & |(\mu-\nu)(\bar P_tf)|=\inf_{\pi\in \scr C(\mu,\nu)}\bigg|\int_{\R^d\times \R^d} \big(\bar P_t f(x)-\bar P_tf(y)\big)\pi(\d x,\d y)\bigg|\\
 &\le  c_2  \W_q(\mu,\nu) \Big[\big(\mu+\nu\big) \Big(\big(\scr M|\nn \bar P_tf|\big)^{\ff q{q-1}}\Big)\Big]^{\ff {q-1}q}+   c_2 t^{-\dd'} \W_1(\mu,\nu),\\
 & \ t\in (0,T],\ z\in \R^d,\ f\in \D_{z,k}.\end{split}\end{equation}
  Moreover,  by \eqref{01},
 \beg{align*}  & \Big[(\mu+\nu)\Big(\big(\scr M|\nn \bar P_tf|\big)^{\ff q{q-1}}\Big)\Big]^{\ff {q-1}q}
\les  \big(\|\mu\|_{\tt L^{p'}}+\|\nu\|_{\tt L^{p'}}\big)^{\ff{q-1}q}
   \Big\|\big(\scr M {|\nn \bar P_t f |} \big)^{\ff q {q-1}}\Big\|_{\tt L^{\ff{p'}{p'-1}}}^{{\ff{q-1}q}}\\
  &\les  \big(\|\mu\|_{\tt L^{p'}}+\|\nu\|_{\tt L^{p'}}\big)^{\ff{q-1}q} \big\|{\nn \bar P_t  }\big\|_{\tt L^{\ff k{k-1}}\to \tt L^{\ff{p'q}{(p'-1)(q-1)}}}\\
  &\les  \big(\|\mu\|_{\tt L^{p'}}+\|\nu\|_{\tt L^{p'}}\big)^{\ff{q-1}q} t^{-\dd'},\ \
   \ t\in (0,T], z\in \R^d,\ f\in \D_{z,k}.\end{align*}
 Combining this with \eqref{02} and $\W_1\le \W_q$, we find a constant $c_2\in (0,\infty)$ such that
   $$\sup_{z\in \R^d,\ f\in \D_{z,k}}  |(\mu-\nu)(\bar P_tf)|\le c_2 \big(\|\mu\|_{\tt L^{p'}}+\|\nu\|_{\tt L^{p'}}\big)^{\ff{q-1}q} t^{-\dd'}\W_q(\mu,\nu),\ \ t\in (0,T], \ \mu,\nu\in \tt{\scr P_{p}}.$$
   Because of  \eqref{dd'},   we  derive \eqref{R1} by repeating the proof of  \eqref{R} with this estimate in place of  \eqref{89}.    When $d=1$ and $\tau>0$ we may take $p'=q=1$, so that
   \eqref{R1'} follows from \eqref{R1}. \end{proof}

 \section{Estimate on the relative entropy}

 \beg{thm}\label{T2.1} Assume {\bf (A)}  and     $\eqref{P3'}$  for some $\tau\in [0, \infty)$,    let $p\in (d,k]\cap [1,\infty)$. Then the following assertions hold.
 \beg{enumerate} \item[$(1)$] For any   $p'\in [1,p]\cap (\ff d{1+2\tau}\lor \ff{dk }{d+k}, p]$, there   exists a constant $c\in (0,\infty)$ such that
\beq\label{ET}\beg{split} \Ent(P_t^*\mu|P_t^*\nu)&\le   \|\mu-\nu\|_{\tt L^{p'}}^2\exp\Big[c + c t( \|\mu\|_{\tt L^{p'}}\land \|\nu\|_{\tt L^{p'}})^{\ff{2p'}{(1+2\tau)p'-d}}\Big]\\
&\quad + \ff {c }t \W_2(\mu,\nu)^2 ,\ \ \ t\in (0,T],\ \mu,\nu\in \tt{\scr P}_p.\end{split} \end{equation}
   In particular, if $d=1$ and $\tau>0$, then there exists a constant $c>0$ such that
 \beq\label{ET'} \Ent(P_t^*\mu|P_t^*\nu)\le  \ff {c }t \W_2(\mu,\nu)^2 + c  \|\mu-\nu\|_{\tt L^{1}}^2,\ \   t\in (0,T],\ \mu,\nu\in \tt{\scr P}_p.\end{equation}
 \item[$(2)$]  For any $p'\in [1,p]\cap (\ff d{2\tau}\lor \ff{dk}{d+k}, k]$, and   $q$  in \eqref{PQ} with $q>\ff{d(p'-1)}{2p'\tau-d},$
 there exists a constant $c>0$ such that
 \beq\label{ET1} \beg{split} \Ent(P_t^*\mu|P_t^*\nu)&\le  \W_q(\mu,\nu)^2
   \exp\Big[c + c t( \|\mu\|_{\tt L^{p'}}\land \|\nu\|_{\tt L^{p'}})^{\ff{2p'}{(1+2\tau)p'-d}}\Big]\\
   &\quad + \ff {c }t \W_2(\mu,\nu)^2 ,\ \ \ t\in (0,T],\ \mu,\nu\in \tt{\scr P}_p.\end{split} \end{equation}
    If in particular $d=1$ and $\tau>\ff 1 2$, then there exists a constant $c>0$ such that
  \beq\label{ET1'} \Ent(P_t^*\mu|P_t^*\nu)\le  \ff {c }t \W_2(\mu,\nu)^2,\ \  \  t\in (0,T],\ \mu,\nu\in \tt{\scr P}_p.\end{equation}
\end{enumerate}
 \end{thm}

 \beg{proof}    (1) For any $\mu\in \tt{\scr P}_p$ and $x\in\R^d$, denote
 $$\rr_t^\mu(x):= \ff{\d P_t^*\mu}{\d x}$$ and
 let
 $X_t^{\mu,x}$ solve the SDE
 $$\d X_t^{\mu,x}= b_t(X_t^{\mu,x}, \rr_t^\mu(X_t^{\mu,x}),\rr_t^\mu) \d t+ \si_t(X_t^{\mu,x})\d W_t,\ \ t\in [0,T],\ X_0^{\mu,x}=x,$$
 which is well-posed under the assumption {\bf (A)}, according to \cite[Proposition  5.1]{HRW25}.
 Let $P_t^\mu f(x):= \E[f(X_t^{\mu,x})]$ for $f\in \B_b(\R^d)$, and
\beq\label{OH} P_t^{\mu,x}:= \L_{X_t^{\mu,x}},\ \ \ P_t^{\mu,\nu}:= \int_{\R^d} P_t^{\mu,x} \nu(\d x),\ \ t\in [0,T],\ \nu\in \scr P.\end{equation}
 We have
\beq\label{PM} P_t^*\mu= P_t^{\mu,\mu}=  \int_{\R^d} P_t^{\mu,x} \mu(\d x),\ \ \ t\in [0,T],\ \mu\in \tt{\scr P}_p.\end{equation}
By \cite[Lemma 2.1]{RW26}, for any $\mu,\nu\in \tt{\scr P}_p$ and $n>1$,
\beq\label{Y0} \Ent(P_t^{\mu,x}|P_t^{\nu, y}) \le n \Ent(P_t^{\mu,x}|P_t^{{\nu, x}})+ (n-1) \log \int_{\R^d} \bigg(\ff{\d P_t^{\nu,x}}{\d P_t^{\nu,y}}\bigg)^{\ff { n }{n-1}} \d P_t^{\nu, y},\ \ t\in (0,T].\end{equation}
According to \cite[Theorem 2.2]{R23}, which applies directly under {\bf (A)} with $l=2$, but the case for $l\ge 3$ can be deduced by induction as in the proof of \cite[Proposition 5.2]{HRW25},
we find constants $c_1>0$ and $n>1$ such that
$$(P_t^{\nu}  f)^n(x)\le \big(P_t^\nu f^n(y) \big) \e^{\ff{c_1|x-y|^2}{t}},\ \ t\in (0,T],\ x,y\in \R^d,\ f\ge 0.$$
By \cite[Theorem 1.4.2(1)]{W13}, this is equivalent to
\beq\label{Y1} ( n-1)\log \int_{\R^d} \bigg(\ff{\d P_t^{\nu,x}}{\d P_t^{\nu,y}}\bigg)^{\ff {n}{n-1}} \d P_t^{\nu, y} \le \ff{c_1|x-y|^2}{t},\ \ x,y\in \R^d,\ t\in (0,T].\end{equation}
To estimate the other term in the upper bound of \eqref{Y0}, we apply Girsanov's transform. For fixed $t\in (0,T],$ let
\beg{align*} &R_t:=\e^{\int_0^t\<\xi_s,\d W_s\>-\ff 1 2\int_0^t |\xi_s|^2\d s},\\
&\xi_s:= (\si_s^*a_s^{-1})(X_s^{\nu,x})\big\{b_s(X_s^{\nu,x},\rr_s^\mu(X_s^{\nu,x}),\rr_s^\mu) -  b_s(X_s^{\nu,x},\rr_s^\nu(X_s^{\nu,x}),\rr_s^\nu)\big\},\ \ s\in [0,t].\end{align*}
By \eqref{P1} and \eqref{P3'}, we find a constant $c_2>0$ such that
\beq\label{YY} |\xi_s|^2\le c_2s^{2\tau}\big(|\rr_s^\mu-\rr_s^\nu|^2(X_s^{\nu,x}) +  s^{-\ff d k} \|\rr_s^\mu-\rr_s^\nu\|_{\tt L^k}^2\big).\end{equation}
Since  $p'>\ff d{1+2\tau}$ implies
$$\ff{2k}{k-d}< \ff{2p'k}{(d(k-p')-2p'k\tau)^+},$$
we find a constant $r\in (\ff{2k}{k-d}, \ff{2p'k}{(d(k-p')-2p'k\tau)^+})$ so that
   $\ff 2 r+\ff d k<1$ and
   $$\dd:= r\Big(\ff{d(k-p')}{2p'k} -\tau\Big)<1.$$
Then by   \eqref{R},  we find a constant  $c_3>0$ such that
\beq\label{C9}  \beg{split}&   \int_0^t \|{s^{\tau}(\rr^\mu_s-\rr^\nu_s)}\|_{\tt L^k}^r\d s\\
&\le \Big(\|\mu-\nu\|_{\tt L^{p'}} \e^{c+c t( \|\mu\|_{\tt L^{p'}}\land \|\nu\|_{\tt L^{p'}})^{\ff{2p'}{(1+2\tau)p'-d}}}\Big)^r
\int_0^t s^{-\dd}\d s\\
&\le  \|\mu-\nu\|_{\tt L^{p'}}^r \e^{c_3+c_3 t( \|\mu\|_{\tt L^{p'}}\land \|\nu\|_{\tt L^{p'}})^{\ff{2p'}{(1+2\tau)p'-d}}},\ \ t\in (0,T].\end{split} \end{equation}
Moreover, $p'>\ff d{1+2\tau}$ implies   $$  \ff{d(k-p')}{ p'k} +\ff  d k-2\tau<1,$$ so  by \eqref{R} we find a constant $c_4>0$ such that
\beq\label{YY1}\int_0^t s^{2\tau -\ff d k} \|\rr_s^\mu-\rr_s^\nu\|_{\tt L^k}^2\d s \le c_4\|\mu-\nu\|_{\tt L^{p'}}^2
 \e^{c_4+c_4 t( \|\mu\|_{\tt L^{p'}}\land \|\nu\|_{\tt L^{p'}})^{\ff{2p'}{(1+2\tau)p'-d}}}. \end{equation}
Combining this with \eqref{YY}, \eqref{C9}  and applying
 Khasminski's estimate, see Theorem \ref{KH} in the next section,  we derive
 $$\E[\e^{\ll \int_0^t |\xi_s|^2\d s}]<\infty,\ \ \ll\in (0,\infty).$$
 So, by Girsanov's theorem,
 $$\tt W_s:= W_s-\int_0^s\xi_r\d r,\ \ s\in [0,t]$$
 is an $m$-dimensional Brownian motion under the probability $\d\Q:= R_t\d \P$. Noting that $X_s^{\nu,x}$ solves the SDE
 $$ \d X_s^{\nu,x}= b_s(X_s^{\nu,x}, \rr_s^\mu(X_s^{\nu,x}), {\rr_s^\mu}) \d s+ \si_s(X_s^{\nu,x})\d \tt W_s,\ \ s\in [0,t],\ X_0^{\nu,x}=x,$$
 by the weak uniqueness, and applying Young's inequality, for any $1<f\in \B_b(\R^d)$, we have
 \beg{align*}P_t^\mu\log f(x) &= \E[R_t \log f(X_t^{\nu,x})] \le \log \E[ f(X_t^{\nu,x})]+\E[R_t\log R_t] \\
&= \log P_t^\nu f(x) + \E[R_t\log R_t].\end{align*}
 Thus, by  \eqref{YY} and  Girsanov's theorem, we obtain
 \beg{align*}&\Ent(P_t^{\mu,x}|P_t^{\nu,x})\le \E[R_t\log R_t]= \ff 1 2 \E_{\Q}  \int_0^t |\xi_s|^2{\d s}\\
 &\le c_2 \E_{\Q} \int_0^t s^{2\tau}  \Big(|\rr_s^\mu-\rr_s^\nu|^2(X_s^{\nu,x}) + s^{-\ff d k}\|\rr_s^\mu-\rr_s^\nu\|_{\tt L^k}^2\Big)\d s.
 \end{align*}
Noting that $\ff 2 r+\ff d k<1$,  by Krylov's estimate, see \cite[Theorem 1.2.3(2)]{RW25}, under {\bf (A)} there exists a constant     $c_3 >0$   such that
\beq\label{KR}  	\beg{split}& \E_{\Q} \int_0^t s^{2\tau}   |\rr_s^\mu-\rr_s^\nu|^2(X_s^{\nu,x})  \d s\leq c_3
	\big\|(\cdot)^{2\tau}(\rr_\cdot ^\mu-\rr_\cdot^\nu)^2\big\|_{\tt L^{k/2}_{r/2}(0,t)}\\
	&\le c_3\bigg(\int_0^t \big(s^{\tau}\|\rr_s ^\mu-\rr_s^\nu\|_{\tt L^{ k }}\big)^r \d s \bigg)^{\ff 2 r},\ \ t\in (0,T].
	\end{split}\end{equation}
 Combining this with \eqref{C9}  {and} \eqref{YY1}, we find   constants $c_4>0$ such that
$$\Ent(P_t^{\mu,x}|P_t^{\nu,x})\le   c_4  \|\mu-{\nu\|_{\tt L^{p'}}^2 \e^{c_4t( \|\mu\|_{\tt L^{p'} }\land \|\nu\|_{\tt L^ {p'}})^{\ff{2p'}{(1+2\tau)p'-d}}}},\ \ t\in (0,T].$$
 This together with  \eqref{Y0} and \eqref{Y1} leads to the following estimate for some constant  $c>0$:
 \beq\label{*DD} \Ent(P_t^{\mu,x}|P_t^{\nu, y}) \le \ff{c|x-y|^2}t+ \|\mu-\nu\|_{\tt L^{p'}}^2 \e^{c+c t( \|\mu\|_{\tt L^{p'}}\land \|\nu\|_{\tt L^{p'}})^{\ff{2p'}{(1+2\tau)p'-d}}},\ \ t\in (0,T].\end{equation}
 Equivalently,
\beq\label{PI} P_t^\mu\log f(x)\le \log P_t^\nu f(y)+ \ff{c|x-y|^2}t+ \|\mu-\nu\|_{\tt L^{p'}}^2 \e^{c+c t( \|\mu\|_{\tt L^{p'}}\land \|\nu\|_{\tt L^{p'}})^{\ff{2p'}{(1+2\tau)p'-d}}}.\end{equation}
 By taking integral with respect to the $\W_2$-optimal coupling of $\mu$ and $\nu$,
 and applying   \eqref{OH} and Jensen's inequality
 $$\int_{\R^d} \big(\log P_t^\nu f(y)\big)\nu(\d y)\le \log \int_{\R^d} \big( P_t^\nu f(y)\big)\nu(\d y)=\log \int_{\R^d} f\d P_t^* \nu,$$
   we derive
\beg{align*}&\int_{\R^d} (\log f)\d P_t^{*}\mu = \int_{\R^d} \big(P_t^\mu \log f(x) \big)\mu(\d x)\\
&\le \int_{\R^d} \big(\log P_t^\nu f(y)\big)\nu(\d y)  + \ff {c}t \W_2(\mu,\nu) + {\|\mu-\nu\|_{\tt L^{p'}}^2\e^{c+c t( \|\mu\|_{\tt L^{p'}}\land \|\nu\|_{\tt L^{p'}})^{\ff{2p'}{(1+2\tau)p'-d}}}}\\
&\le \log \int_{\R^d} f\d P_t^* \nu  + \ff {c}t \W_2(\mu,\nu) + \|\mu-\nu\|_{\tt L^{p'}}^2 \e^{c+ ct( \|\mu\|_{\tt L^{p'}}\land \|\nu\|_{\tt L^{p'}})^{\ff{2p'}{(1+2\tau)p'-d}}}.\end{align*}
This implies \eqref{ET}.  When  $d=1$ and $\tau>0$, we may take $p'=1$ so that \eqref{ET} implies \eqref{ET'}.

(2)  Since   $p'> \ff d {2\tau}$ and $q> \ff{d(p'-1)}{2 p'\tau-d}$, we have
\beq\label{TQ}\ff{2k}{k-d}< \Theta:=\ff 1 {[\ff 1 2 +\ff d 2(\ff 1 {p'} +\ff{p'-1}{p'q} -\ff 1 k)-\tau]^+}.\end{equation}
Taking $r\in (\ff{2k}{k-d},\Theta)$, we have  $\ff 2 r+\ff d k<1$ and
\beq\label{DE} \dd:=r\bigg(\ff 1 2 +\ff d 2\Big(\ff 1 {p'} -\ff 1 k +\ff{p'-1}{p'q}\Big)-\tau\bigg)<1.\end{equation}
Then by \eqref{R1}, we find a constant $c_1\in (0,\infty)$ such that
\beq\label{C9'} \beg{split}&\int_0^t \|{s^{\tau}(\rr^\mu_s-\rr^\nu_s)}\|_{\tt L^k}^r\d s\\
&\le \W_q(\mu,\nu)^r(\|\mu\|_{\tt L^{p'}}+\|\nu\|_{\tt L^{p'}})^{\ff{r(q-1)}q}  \e^{c_1+c_1 t( \|\mu\|_{\tt L^{p'}}\land \|\nu\|_{\tt L^{p'}})^{\ff{2p'}{(1+2\tau)p'-d}}},\ \
t\in (0,T].\end{split}\end{equation}
Combining this with   \eqref{KR},  we find  $c_2\in (0,\infty)$ such that
	\begin{equation}\label{10'} \beg{split}
		&\E_\Q \int_0^t s^{2\tau}  |\rr_s^\mu-\rr_s^\nu|^2(X_s^{\nu,x}) \d s\\
		 &\le \W_q(\mu,\nu)^2(\|\mu\|_{\tt L^{p'}}+\|\nu\|_{\tt L^{p'}})^{\ff{2(q-1)}q}  \e^{c_2+c_2 t( \|\mu\|_{\tt L^{p'}}\land \|\nu\|_{\tt L^{p'}})^{\ff{2p'}{(1+2\tau)p'-d}}},\ \ \ \  t\in (0,T].  \end{split}\end{equation}
	Moreover,  $q> \ff{d(p'-1)}{2 p'\tau-d}$ implies
			$$\tt \dd:=  1+d\Big(\ff 1 {p'} +\ff{p'-1}{p' q}  \Big)-2\tau <1.$$ Then, by \eqref{R1} we find   constants     $c_3,c_4\in (0,\infty)$   such that
	\beq
	\beg{split}\label{10}
	&\int_0^t s^{2\tau-\ff d k}   \|\rr_s^\mu-\rr_s^\nu\|_{\tt L^k}^2\d s\\
	&	 \le    \W _q(\mu,\nu)^2 \e^{c_3+c_3 t( \|\mu\|_{\tt L^{p'}}\land \|\nu\|_{\tt L^{p'}})^{\ff{2p'}{(1+2\tau)p'-d}}}
	\int_0^t    s^{ -\tt\dd  }\d s
	\\&\le  \W _q(\mu,\nu)^2 \e^{c_4+c_3 t( \|\mu\|_{\tt L^{p'}}\land \|\nu\|_{\tt L^{p'}})^{\ff{2p'}{(1+2\tau)p'-d}}},\ \ \ t\in (0,T].
	\end{split}
\end{equation}
By repeating the proof of \eqref{ET} using  {\eqref{10'} and \eqref{10} in place of \eqref{C9} and \eqref{YY1}},   we derive   \eqref{ET1}.   When $d=1$ and $\tau>\ff 1 2$,
we may take $p'=1$ and $q=1$ so that  \eqref{ET1'} follows from \eqref{ET1}.
\end{proof}

\section{Estimate on the Renyi entropy}
It is clear that  $\Ent_\aa$ is {increasing} in $\aa$ and
      $$\Ent(\mu|\nu)=\lim_{\aa\downarrow 0} \Ent_\aa(\mu|\nu).$$
      We intend to estimate $\Ent_\aa(P_t^*\mu|P_t^*\nu)$ by using suitable distances  of the initial distributions $\mu$ and $\nu$.
   For any constant $c>0$, let
   $$ \W_{ c,\e}  (\mu,\nu):=\inf_{\pi\in \C(\mu,\nu)} \log \int_{\R^d\times\R^d} \e^{c|x-y|^2}{\pi(\d x,\d y)}. $$

\beg{thm}\label{T3.1} Assume {\bf (A)} and let $p\in (d,k]\cap [1,\infty).$   If $\eqref{P3'}$ holds for some $\tau\in [0, \infty)$, then  the following assertions hold.
 \beg{enumerate} \item[$(1)$] Let  $p'\in [1,p]\cap (\ff d{1+2\tau}\lor \ff{dk}{d+k}, p].$ Then there    exist constants $\aa,c\in (0,\infty)$ and a map
 $\bb:   (\ff {2k}{k-d}, \infty)\to (0,\infty)$ such that for any $r\in (\ff {2k}{k-d},  \infty), t\in (0,T]$ and $\mu,\nu\in \tt{\scr P}_p,$
 \beg{align*} &\Ent_\aa(P_t^*\mu|P_t^*\nu)\le \ff{1}{\alpha} \W_{\ff c {2t},\e}  (\mu,\nu)\\
 &+  \bb(r)\exp\bigg[ ct(\|\mu\|_{\tt L^{p'}}\land \|\nu\|_{\tt L^{p'}})^{\ff{2p'}{(1+2\tau)p'-d}}\bigg]  
    \Big( \|\mu-\nu\|_{\tt L^{p'}}^r +   \|\mu-\nu\|_{\tt L^{p'}}^2\Big).
  \end{align*} If  $d=1$ and $\tau > 0,$ then this estimate holds for $p'=q=1$ such that for some different  $\bb: (\ff {2k}{k-d},  \infty)\to (0,\infty)$,
  $$\Ent_\aa(P_t^*\mu|P_t^*\nu)\le \ff{1}{\alpha} \W_{\ff c {2t},\e} (\mu,\nu)+
   \bb(r)  \Big( \|\mu-\nu\|_{\tt L^{1}}^r +    \|\mu-\nu\|_{\tt L^{1}}^2\Big).$$
   \item[$(2)$] For any $p'\in [1,p]\cap (\ff d{2\tau}\lor\ff{dk}{d+k},  k]$  and  $q$ in $\eqref{PQ}$  with $q> \ff{d(p'-1)}{2p'\tau-d},$ there    exist $\aa,c\in (0,\infty)$ and
 $\bb:   (\ff {2k}{k-d},  \infty) \to (0,\infty)$ such that for any $r  \in (\ff {2k}{k-d},  \infty), t\in (0,T]$ and $\mu,\nu\in \tt{\scr P}_p,$
 \beg{align*} &\Ent_\aa(P_t^*\mu|P_t^*\nu)\le \ff{1}{\alpha} \W_{\ff c {2t},\e}  (\mu,\nu)\\
 &+ \bb(r) \exp\bigg[ c  t(\|\mu\|_{\tt L^{p'}}\land \|\nu\|_{\tt L^{p'}})^{\ff{2p'}{(1+2\tau)p'-d}}\bigg] 
    \Big( \W_q(\mu,\nu)^r +   \W_q(\mu,\nu)^2\Big).
  \end{align*}  If  $d=1$ and $\tau > \ff 1 2,  $ then this estimate holds for $p'=q=1$ such that for some different  $\bb: \in (\ff {2k}{k-d},  \infty)\to (0,\infty)$,
  \beg{align*} \Ent_\aa(P_t^*\mu|P_t^*\nu)\le \ff{1}{\alpha} \W_{\ff c {2t},\e}  (\mu,\nu)
   + \bb(r)  \Big(  \W_1(\mu,\nu)^r +   \W_1(\mu,\nu)^2\Big).
  \end{align*}
\end{enumerate}
 \end{thm}

\beg{proof}  Let $\Theta>\ff{2k}{k-d}$ be in \eqref{TQ}. Noting that
$$\|\cdot\|_{\tt L^{p'}}^{r_1}+\|\cdot\|_{\tt L^{p'}}^{2}\ge \|\cdot\|_{\tt L^{p'}}^{r_2},\ \ \ r_1\ge r_2\ge 2,$$
it suffices to find $\bb:   (\ff{2k}{k-d}, \Theta)\to (0,\infty)$ such that the desired estimates hold for $r\in (\ff{2k}{k-d}, \Theta).$ In this case,  we have already proved \eqref{C9}, \eqref{YY1},  \eqref{C9'} and \eqref{10}  in the proof of Theorem \ref{T2.1}.

(1) By \cite[Theorem 2.2]{R23}, under {\bf (A)}, there exist  constants    $\aa_1 \in (1,\infty) $ and $c_1\in (0,\infty)$  such that
$$(P_t^\mu f(y))^{n} \e^{-\ff{c_1|x-y|^2}t} \le P_t^\mu f^{n}(x),\ \ t\in (0,T],\ x,y\in \R^d,\ \ n\ge \aa_1.$$
Combining this with
  {Girsanov}'s theorem used in the the proof of Theorem \ref{T2.1},   and applying the Schwarz inequality,  we obtain that
for any $0\le f\in \B_b(\R^d)$,
\beq\label{Z0} \beg{split} &(P_t^\mu f(y))^{2n} \e^{\ff{-2c_1|x-y|^2}t} \le \big(P_t^\mu f^{n}(x)\big)^{2}=\Big( \E\big[f^{n}(X_t^{\mu,x})\big]\Big)^{2}\\
&= \Big( \E\big[R_t f^{n}(X_t^{\nu,x})\big]\Big)^{2}\le \Big( \E\big[  f^{2n}(X_t^{\nu,x})\big]\Big) \Big({\E\big[R_t^{2}}\big]\Big).\end{split}\end{equation}
By $\eqref{P3'}$ and the refined  Khasminskii estimate in Theorem \ref{KH}, we find a constants $c_2 >0$ such that
\beq\label{z0'} \beg{split} & { \E\big[R_t^2 ] }\le \exp\bigg[c_2\int_0^t s^{r\tau}\|\rr_s^\mu-\rr_s^\nu\|_{\tt L^k}^r\d s  + c_2\bigg(\int_0^t s^{r\tau}\|\rr_s^\mu-\rr_s^\nu\|_{\tt L^k}^r\d s\bigg)^{\ff 2 r}\\
&\qquad\qquad \qquad\qquad \qquad\qquad  \qquad\qquad+ c_2\int_0^t s^{2\tau-\ff d k}\|\rr_s ^\mu-\rr_s^\nu\|_{\tt L^k}^2\d s\bigg].\end{split}\end{equation}
Combining this with \eqref{C9}, \eqref{YY1} and \eqref{Z0},  we find constants $c>0$ and $\bb(r)\in (0,\infty)$  such that
$$\big(P_t^\mu f(y)\big)^{2n} \e^{-\ff{c(2n-1)|x-y|^2}{2   t}}   \le \big(P_t^\nu f^{2n}(x)\big) \e^{\bb(r) H_{t,\mu,\nu}^c},$$
where
$$H_{t,\mu,\nu,c}:= \big(\|\mu-\nu\|_{\tt L^{p'}}^2+ \|\mu-\nu\|_{\tt L^{p'}}^r\big) \e^{c t (\|\mu\|_{\tt L^{p'}}\land \|\nu\|_{\tt L^{p'}})^{\ff{2p'}{(1+2\tau)p'-d}}}.$$
Integrating both sides with respect to the optimal coupling $\pi(\d y,\d x)$ for $\mu$ and $\nu$ reaching $\W_{\ff c {2t},\e}(\mu,\nu)$, and applying H\"older's inequality, we derive
\beg{align*} &\bigg(\int_{\R^d} f \d P_t^*\mu\bigg)^{2n} = \bigg(\int_{\R^d\times\R^d} (P_t^\mu f)(y)  \pi(\d y,\d x)  \bigg)^{2n}  \\
&\le \bigg(\int_{\R^d\times\R^d } (P_t^\mu f(y))^{2n} \e^{-\ff{c(2n-1)|x-y|^2}{2  t}}  \pi(\d y,\d x)\bigg)  \bigg(\int_{\R^d\times\R^d} \e^{\ff {c|x-y|^2} {2t}} {\pi(\d y,\d x)}\bigg)^{2n-1}\\
&\le     \exp \Big[ \log \int_{\R^d\times\R^d } (P_t^\nu f^{2n} (x)  \e^{\bb(r) H_{t,\mu,\nu}^c} \pi(\d y,\d x) \Big]\exp  \Big[(2n-1)  \W_{\ff c {2t},\e}(\mu,\nu)\Big]
\\
&\le   \exp \Big[     \bb(r) H_{t,\mu,\nu}^c +  \log \int_{\R^d\times\R^d } (P_t^\nu f^{2n} (x)   \pi(\d y,\d x) \Big]  \exp  \Big[(2n-1)  \W_{\ff c {2t},\e}(\mu,\nu)\Big] \\
&\le \exp  \Big[(2n-1) \W_{\ff c{2t},\e}(\mu,\nu) +\bb(r)  H_{t,\mu,\nu}^c\Big]  \int_{\R^d\times\R^d } (P_t^\nu f^{2n} (x)   \pi(\d y,\d x) \\
&=  \exp\Big[(2n-1) \W_{\ff c{2t},\e}(\mu,\nu) +\bb(r)  H_{t,\mu,\nu}^c \Big]   \int_{\R^d  } f^{2n} \d P_t^*\nu ,\ \ n\ge \aa_1.\end{align*}
By \cite[Theorem 1.4.2]{W13}, this implies  the   desired estimate on $\Ent_\aa(P_t^*\mu|P_t^*\nu)$ for $\aa= \ff 1 {2n-1}.$

(2) The proof of the second assertion  is completely similar by using \eqref{C9'}    and \eqref{10}  in place of \eqref{C9} and \eqref{YY1}.
\end{proof}

   \section{Refined Khasminskii estimate }

The  Khasminskii estimate, which goes back to \cite{K},  is a power tool in the study of singular SDEs.

 Let $T\in (0,\infty)$.   Consider the following SDE on $\R^d$:
\beq\label{SDE} \d X_t= b_t(X_t)\d t+\si_t(X_t)\d W_t,\ \ \ t\in [0,T],\end{equation}
where $W_t$ is an $m$-dimensional Brownian motion on a complete filtered probability space $(\OO, \F, \{\F_t\}_{t\in [0,T]},\P),$ and
$$b: [0,T]\times \R^d\to \R^d,\ \ \ \si: [0,T]\times\R^d\to \R^{d\otimes m}$$
are measurable satisfying the following conditions.

    \beg{enumerate}
\item[{\bf (B)}]       There exist a constant $K\in (0,\infty),l\in \mathbb N$ and  $\{(p_i,q_i)\}_{1\le i\le l}\subset \scr K$ such that $\si$ and $b$ satisfy the following conditions on $[0,T]\times\R^d$.
\item[$(1)$]     $a:= \si\si^*$ is invertible with $\|a\|_\infty+\|a^{-1}\|_\infty\le K$, where $\si^*$ is the transposition of $\si$,   and
$$   \zeta(\vv):=\sup_{|x-y|\le \vv, t\in [0,T]} \|a_t(x)-a_t(y)\| \downarrow 0\ \text{as}\ \vv\downarrow 0.$$
\item[$(2)$]  $b=\sum_{{i=1}}^l b^{(i)}$, $b^{(1)}$ is locally bounded and
 \beq\label{LPS}  \|\nn b^{(1)}\|_\infty+ \sum_{i=2}^l \|b^{(i)}\|_{\tt L_{q_i}^{p_i}(T)}  \le K.\end{equation}
 \end{enumerate}

By the Khasminskii  estimate,  see  \cite[Theorem 1.2.3(2), Theorem 1.2.4]{RW25},  under this condition with $l=2$, for any $(p,q)\in \scr K$, there exist constant $c>0$ and $\kk >2$ depending only on $d,K,T, p,q$ and
$\zeta$,
 such that
\beq\label{*1} \E\big(\e^{\int_s^t f(X_r)^2\d r}\big|\F_s\big) \le \e^{c  +c  \|f\|_{\tt L_q^p(s,t)}^\kk}.\end{equation}
The earlier versions of this type estimate are given for $f\in L_q^p(s,t):= L^q([s,t]\to L^p(\R^d))$. The version with $f\in \tt L_q^p(s,t)$ is first proved in \cite{XXZZ}  under {\bf (B)} with $l=2$ and $b^{(1)}=0$,
which is then extended in \cite{YZ} to $b^{(1)}\ne 0$ and $l=2$.  We present below a refined version with $\kk=q$.

\beg{thm}\label{KH} Assume $ {\bf (B)}$. Then for any $(p,q)\in \scr K$, there exists a constant $c\in (0,\infty)$ depending only on $d,K,T, p,q$ and
$\zeta$ such that any solution to $\eqref{SDE}$ satisfies
\beq\label{ES1}\E\big(\e^{\int_s^t f_r(X_r)^2\d r} \big|\F_s\big)\le \beg{cases}
\e^{c\|f\|_{\tt L_q^p(s,t)}^2},\ &\text{if}\ \|f\|_{\tt L_q^p(s,t)}\le 1,\\
\e^{c \int_s^t \|f_r\|_{\tt L^p}^q\d r},\ &\text{otherwise} \end{cases}\end{equation}
for any $0\le s<t\le T$ and $f\in \tt L_q^p(s,t).$
Consequently,
\beq\label{ES2} \E\big(\e^{\int_s^t f_r(X_r)^2\d r} \big|\F_s\big)\le \e^{c\|f\|_{\tt L_q^p(s,t)}^2
+ c \int_s^t \|f_r\|_{\tt L^p}^q\d r},\ \ 0\le s<t\le T,\ f\in \tt L_q^p(s,t).\end{equation}
\end{thm}

\beg{proof} All constants below depend only on $d,K,T, p,q$ and
$\zeta$.

(a) Let {$l=2$}. In this case, the assumption {\bf (B)} coincides with $(A^{1.1})$ in
\cite{RW25}.
  By Jensen's inequality and \eqref{*1} for $\ll f$ in place of $f$,  we find a constant $c_1>0$ and $\kk >2$  such that
\beg{align*}  \Big[\E\big(\e^{\int_s^t f(X_r)^2\d r}\big|\F_s\big)\Big]^{\ll^2} \le
\E\big(\e^{\int_s^t (\ll f)(X_r)^2\d r}\big|\F_s\big)
 \le \e^{c_1 +c_1 \ll^\kk\|f\|_{\tt L_q^p(s,t)}^\kk},\ \ \ll\ge 1,\end{align*}
 so that
 $$\E\big(\e^{\int_s^t f(X_r)^2\d r}\big|\F_s\big)\le \e^{c_1\ll^{-2}  +c_1 \ll^{\kk-2}\|f\|_{\tt L_q^p(s,t)}^\kk},
 \ \ \ll\ge 1.$$
 By taking    $\ll=   (1\land \|f\|_{\tt L_q^p(s,t)})^{-1},$ we derive
\beq\label{*3} \E\big(\e^{\int_s^t f(X_r)^2\d r}\big|\F_s\big)\le
 \e^{2c_1\|f\|_{\tt L_q^p(s,t)}^2+ c_1 \|f\|_{\tt L_q^p(s,t)}^\kk}.\end{equation}
Consequently, when $\|f\|_{\tt L_q^p(s,t)}\le 1$, \eqref{ES1} holds for $c=2c_1.$  

 Now, let $\|f\|_{\tt L_q^p(s,t)}>1$. We have
\beq\label{A*} A:=\int_s^t \|f_r\|_{\tt L^p}^q\d r \ge\|f\|_{\tt L_q^p(s,t)}^q>1.\end{equation}
 For any $1\le n\in \mathbb N,$ choose $s=t_0<t_1<\cdots<t_n=t$ such that
 $$\|f\|_{\tt L_q^p(t_i,t_{i+1})}^q \le \int_{t_i}^{t_{i+1}}\|f_r\|_{\tt L^p}^q\d r = \ff 1 n A,\ \ \ 0\le i\le n-1.$$
Combining this with \eqref{*3} for $(t_i, t_{i+1})$ in place of $(s,t)$, we derive
 $$\E\Big(\e^{\int_{t_i}^{t_{i+1}} f(X_r)^2\d r}\Big|\F_s\Big)\le  \e^{2 c_1 n^{-\ff 2q} A^{\ff 2 q} + c_1 n^{-\ff \kk q}A^{\ff\kk q}},\ \ 0\le i\le n-1.$$
 Hence,
 \beg{align*} &\E\Big(\e^{\int_s^t f(X_r)^2\d r}\Big|\F_s\Big)=
 \E\Big[\e^{\int_s^{t_{n-1}} f(X_r)^2\d r} \E\Big(\e^{\int_{t_{n-1}}^t f(X_r)^2\d r}\Big|\F_{t_{n-1}}\Big) \Big|\F_s\Big]\\
 &\le \e^{2 c_1 n^{-\ff 2q} A^{\ff 2 q} + c_1 n^{-\ff \kk q}A^{\ff\kk q}} \E\big(\e^{\int_s^{t_{n-1}} f(X_r)^2\d r}\big|\F_s\big)\\
 &\le \cdots\le \e^{2 c_1 n^{1-\ff 2q} A^{\ff 2 q} + c_1 n^{1-\ff \kk q}A^{\ff\kk q}}.\end{align*}
 Taking
 $$n= \inf\big\{m\in\mathbb N:\ m\ge A\big\},$$
which satisfies $2A\ge n \ge A$  due to $A>1$,  we find a constant $c >0$ such that
$$\E\big(\e^{\int_s^t f(X_r)^2\d r}\big|\F_s\big)\le \e^{c  A}= \e^{c\int_s^t \|f_r\|_{\tt L_p}^q\d r},$$
so that \eqref{ES1} holds for $\|f\|_{\tt L_q^p(s,t)}>1.$

 (b) Suppose that the estimate holds for $l=l_0$ for some $l_0\in \mathbb N$, it remains to prove
 it for $l=l_0+1.$ This can be done by the Zvoinkin's transform. More precisely, let
 $$L_t:=\ff 1 2 {\rm tr}\big\{\si_t\si_t^*\nn^2\big\}+\sum_{i=1}^{l_0+1}b_t^{i} \cdot\nn,\ \ t\in [0,T].$$
 By \cite[Theorem 1.2.3(3)]{RW25}, see also \cite{YZ}, when $\ll>0$ is large enough, the PDE
 $$(\pp_t +L_t)u_t= \ll u_t- b_t^{(l_0+1)},\ \ \ t\in [0,T],\ u_T=0$$
 has a unique solution such that
 $ \|\nn^2 u\|_{\tt L_{\infty}^{p}(T)}<\infty$ for some $p>d,$ and
 $$\|u\|_\infty+\|\nn u\|_\infty\le \ff 1 2.$$
 So, $\nn\Theta_t(x)$ is H\"older continuous in $x$ uniformly in $t\in [0,T]$ and
 $$\Theta_t(x):= x+u_t(x),\ \ t\in [0,T]$$ are diffeomorhisms with $\|(\nn \Theta)^{-1}\|_\infty<\infty.$
 By It\^o's formula, see  \cite[Theorem 1.2.3(3)]{RW25},
 $Y_t:=\Theta_t(X_t)$ satisfies
 $$\d Y_t = \tt b_t(Y_t)\d t+\tt \si_t(Y_t)\d W_t,$$
 where
 $$\tt \si_t:= \big\{(\nn \Theta_t)^*\si_t\big\}\circ\Theta_t^{-1},\ \
 \tt b_t:=\Big\{\ll u_t+ \sum_{i=1}^{l_0}  b_t^{(i)} \Big\}\circ\Theta_t^{-1}$$
 satisfy \cite[$(A^{1,1})$]{RW25}  in stead of $\si$ and $b$ for $l=l_0.$
 Then the proof is finished by the desired estimate for $l=l_0,$
and  for $(Y_r, f_r\circ\Theta_r^{-1})$ in place of $(X_r, f_r)$.

\end{proof}

%\bibliographystyle{plain}
%\bibliography{zyPHD}

\end{document}